\documentclass{amsart}
\usepackage{graphicx}
\usepackage{amscd}
\usepackage{url}
\usepackage{amssymb}

\newtheorem{theorem}{theorem}[section]
\newtheorem{lemma}[theorem]{Lemma}

\theoremstyle{definition}
\newtheorem{definition}[theorem]{Definition}
\newtheorem{example}[theorem]{Example}

\theoremstyle{remark}
\newtheorem{remark}[theorem]{remark}

\numberwithin{equation}{section}

\newcommand{\ip}[2]{\langle#1,#2\rangle}
\newcommand{\nm}[1]{\|{#1}\|}

\numberwithin{equation}{section}

\DeclareMathOperator{\PFP}{PFP}
\DeclareMathOperator{\FP}{FP}

\newcommand{\C}{\mathbb{C}}
\newcommand{\R}{\mathbb{R}}

\newcommand{\N}{\mathbb{N}}

\newcommand{\Vol}{\operatorname{Vol}}
\renewcommand{\chi}{{\mathbbm{1}}}

\newcommand{\sca}{\mathcal{S}\mathcal{C}}
\newcommand{\calC}{\mathcal C}

\DeclareMathOperator{\supp}{supp}

\newcommand{\fr}{\mathcal{F}}

\newcommand{\co}{\operatorname{co}}

\newcommand{\aff}{\operatorname{aff}}

\DeclareMathOperator{\pfp}{PFP}
\DeclareMathOperator*{\linspan}{span}

\newcommand{\tr}{\text{Tr}}

\newtheorem{proposition}[theorem]{proposition}
\newtheorem{cor}[theorem]{Corollary}
\newtheorem{question}[theorem]{Question}

\begin{document}
\title{Preconditioning techniques in frame theory and probabilistic frames}
\author{Kasso A.~Okoudjou}
\address{
Department of Mathematics $\&$ \\
Norbert Wiener Center\\
University of Maryland\\
College Park, MD 20742 USA}
\email{kasso@math.umd.edu}

\date{\today}

\maketitle

\begin{abstract}
In this chapter we survey two  topics that have recently been investigated in frame theory.   First, we give an overview of the class of scalable frames. These are (finite) frames with the property that each frame vector can be rescaled in such a way that the resulting frames are tight. This process can be thought of  as a preconditioning method for finite frames. In particular, we: (1) describe the  class of  scalable frames; (2) formulate various equivalent characterizations of scalable frames, and relate the scalability problem to the Fritz John ellipsoid theorem. Next, we  discuss  some results on a  probabilistic interpretation of frames. In this setting, we: (4) define probabilistic frames as a generalization of frames and as a subclass of continuous frames; (5) review the properties of certain potential functions whose minimizers are frames with certain optimality properties. 
\end{abstract}

\tableofcontents

\section{Introduction}

This chapter  is devoted to two topics that have been recently investigated within frame theory:  (a) frame preconditioning methods  investigated under the vocable \emph{scalable frames}, and (b) probabilistic methods in frame theory referred to as \emph{probabilistic frames.} Before getting to the details on each of these topics we briefly review some essential facts on frame theory and refer to  \cite{framesbook, heil13, koche1, koche2} for more on frames and their applications. In all that follows we restrict ourselves to (finite)  frames in $\R^N$. 

\subsection{Review on finite frame theory}\label{subsec1.1}

\begin{definition}\label{framedef}
A set $\Phi=\{\varphi_k\}_{k=1}^{M}
\subseteq \R^N$ is a \emph{frame}\index{frame} for $\R^N$ if
$   \exists \, A, B >0 $ such that $\forall x \in \R^N,$ $$
    A\nm{x}{}^{2} \leq \sum_{k=1}^{M} | \ip{x}{\varphi_k}|^{2}
   \leq B \nm{x}{}^{2}.
$$
If, in addition, each $\varphi_k$ is unit-norm, we
say that $\Phi$ is a {\it unit-norm frame}. The set of frames
for $\R^N$ with $M$ elements will be denoted by $\fr(M,N)$, and simply $\fr$ if $M$ and $N$ are fixed. In addition, we shall denote the subset of unit-norm frames by $\fr_u(M,N)$, i.e.,  
$$
\fr_u(M,N) := \left\{\{\varphi_k\}_{i=1}^M\in\fr(M,N) : \|\varphi_k\|_2=1\text{ for }k=1,\ldots,M\right\}.
$$
\end{definition}

We shall investigate frames via the {\it analysis operator}, $L$,  defined by

$$ L:  \R^N \to \R^M: 
  x\mapsto  L x=\{\ip{x}{\varphi_{k}}\}_{k=1}^{M}.$$ The {\it synthesis operator} is the adjoint  $L^*$ of $L$ and is defined by $$L^*: \R^M \to \R^N:   c=(c_k)_{k=1}^M \mapsto L^{*}c=\sum_{k=1}^{M}c_k \varphi_k.$$ It is easily seen that the canonical matrix associated to $L^*$ is the $N\times M$ matrix whose $k^{th}$ column is the frame vector $\varphi_k$. As such we shall abuse notation and denote this matrix by $\Phi$ again. Consequently, the canonical matrix associated with $L$ is simply $\Phi^T$, the transpose of $\Phi$. 

The {\it frame operator}\index{frame!frame operator} $S=L^{*}L$ is given by $$S: \R^N \to \R^N:  x \mapsto Sx=\sum_{k=1}^M\ip{x}{\varphi_k}\varphi_k ,$$ and its matrix will be denoted (again) by $S$ with $$S=\Phi \Phi^{T}.$$  The {\it Gramian (operator)}\index{Gram operator} of the frame  is defined by $$G=L L^{*}=\Phi^T\Phi.$$ In fact, the Gramian in an $M\times M$ matrix whose $(i, j)^{th}$ entry is $\ip{\varphi_{j}}{\varphi_i}.$  

$\Phi=\{\varphi_k\}_{k=1}^{M}\subset \R^M$ is a frame if and only if $S$ is a positive definite  matrix on $\R^M$. In this case, $$\{\tilde{\varphi}_{k}\}_{k=1}^{M}=\{S^{-1}\varphi_{k}\}_{k=1}^{M}$$ is also a frame, called \emph{the canonical dual frame}\index{frame!canonical dual frame}, and, for each $x \in \R^N$, we have
\begin{equation}\label{framerecons}
x = \sum_{k=1}^{M}\ip{x}{\varphi_{k}}\tilde{\varphi}_{k}=\sum_{k=1}^{M}\ip{x}{\tilde{\varphi}_{k}}\varphi_{k}.
\end{equation}

A frame $\Phi$ is a \emph{tight frame}\index{frame!tight frame} if we can choose $A = B$.  In this case the frame operator is simply a multiple of the identity operator. To any frame  $\Phi=\{\varphi_{k}\}_{k=1}^{M} \subset \R^M$ is associated a \emph{canonical tight frame} given by $$\{\varphi_{k}^{\dagger}\}_{k=1}^{M}=\{S^{-1/2}\varphi_{k}\}_{k=1}^{M}\subset \R^N$$ such that for every $x\in \R^N$, 
\begin{equation}\label{tfrecon2}
x=\sum_{k=1}^{M}\ip{x}{\varphi_{k}^{\dagger}}\varphi_{k}^{\dagger}.
\end{equation}

 If   $\Phi$ is a tight frame of unit-norm vectors,  we say that $\Phi$ is  a \emph{finite unit-norm tight frame (FUNTF)}\index{frame!FUNTF}. In this case,   the reconstruction formula~\eqref{framerecons} reduces  to 
\begin{equation}\label{tframerecons}
\forall x \in \R^N, \quad x = \tfrac{N}{M}\sum_{k=1}^{M}\ip{x}{\varphi_{k}}\varphi_{k}.
\end{equation} 
FUNTFs are one of the most fundamental objects in frame theory with several applications.  Other chapters  in this volume will delve more into FUNTFs. 
The reconstruction formulas~\eqref{tframerecons} and~\eqref{tfrecon2}  are  very reminiscent of the expansion of a vector in an orthonormal basis for $\R^d$.  However, due to the redundancy of frames, the coefficients in these reconstruction formulas are not unique. But the ``simplicity'' of these  reconstruction formulas  makes the use  of tight frames very attractive in many applications. This in turn, spurs the need for methods to characterize and construct tight frames.

A major development in this direction is due to Benedetto and Fickus~\cite{bf03} who proved that  for each  $\Phi=\{\varphi_k\}_{k=1}^{M}\subset \R^N$, such that $\|\varphi_k\|=1$ for each $k$, we have 
\begin{equation}\label{potineq}
\text{FP}(\Phi) =\sum_{j=1}^{M} \sum_{k=1}^{M} |\ip{\varphi_{j}}{\varphi_{k}}|^{2}\geq \tfrac{M}{N}\max(M, N),
\end{equation}  where $\text{FP}(\Phi)$ is the \emph{frame potential.} \index{frame!frame potential}
The bound given in ~\eqref{potineq} is the global minimum of $\text{FP}$  and is achieved by orthonormal systems if $M\leq N$, and by  tight frames if $M>N$.
Casazza, Fickus, Kova\v{c}evi\'{c}, Leon and Tremain \cite{CFKLT}  extended this result by removing the condition that the frames are unit norm.  In essence, these results suggest that one may effectively search for FUNTFs by minimizing the frame potential.  In practice, techniques such as the steepest-descent method can be used to find the minimizers of the frame potential  \cite{cfm11, massruiz}. 
For other related results on the frame potential we refer to \cite{fjko, jbok}, and  \cite[Chapter 10]{framesbook}. The frame potential and some of its generalization will be considered in Section~\ref{sec3}.

Construction of FUNTFs has seen a lot of research activities in recent years and as a result a number of construction methods have been offered. Casazza, Fickus, Mixon, Yang and Zhou \cite{cfmyz11}  introduced the \emph{spectral tetrix} method  for constructing FUNTFs.  This method has been extended to construct all tight frames,  \cite{chkwz, fmps11a, lmo12}.  There have also been other new insights in the construction of tight frames, leading to methods rooted in differential and algebraic geometry, see \cite{cfmps13, cast13, stra12, stra11}. Some of these methods will be introduced in some of the other chapters  of this volume. 

\subsection{Scalable frames}\label{subsec1.2}
While these powerful algebro-geometric methods can construct all tight frames with given vector lengths and frame constant, they have not been able to incorporate any extra requirement. To put it simply, constructing application-specific FUNTFs involves extra construction constraints, which usually makes the problem very difficult. However, one could ask if  a (non) tight frame can be (algorithmically) transformed into a tight one. An analogous problem has been investigated for decades in numerical linear algebra. Indeed, preconditioning methods are routinely used to convert large and poorly-conditioned systems of linear equations  $Ax=b$, into better conditioned ones  \cite{ben, forstraus,  nhi02}. For example,  a matrix $A$ is (row/column) scalable if there exit diagonal matrices $D_1, D_2$ with positive diagonal entries such that $D_1A, AD_2,$ or $D_1AD_2$ have  constant row/column sum, \cite{babo95, forstraus, jr09, khka92,  stwitz62}. Matrix scaling is part of general preconditioning schemes in numerical linear algebra \cite{ben, nhi02}. 

One of the goals of these notes  is to  survey recent developments in preconditioning methods in  finite frame theory. In particular, we describe recent developments in answering the following question:

\begin{question}\label{q1}
 Is it possible to (algorithmically) transform a (non-tight) frame into  a tight frame? 
\end{question}

In Section~\ref{sec2}, we  outline  a  convex geometric approach that has been proposed to  answer this question.   For example, we consider the case of  solving Question~\ref{q1} using some classes of scaling matrices to transform a non-tight frame into a tight one.  More specifically, we give an overview of recent results addressing  the following  problem.

\begin{question}\label{q2}
 Is it possible to rescale the norms of the vectors in a (non-tight) frame to obtain a tight frame? 
\end{question}

Frames that answer positively Question~\ref{q2} are termed  \emph{scalable frames}\index{frame!scalable frames}, and  were  characterized in \cite{kopt12}, see also \cite{kop13a}. This characterization is operator-theoretical and solved the problem in both the finite and the infinite dimensional settings. More precisely,  in the finite dimensional setting, the main result of \cite{kopt12} characterizes the set of non scalable frames and gives a simple geometric condition for a frame to be non scalable in dimensions $2$ and $3$, see Section~\ref{subsec2.2}. Other characterizations of scalable frames  using the properties of the so-called diagram vector\index{diagram vector} (\cite{hklw}) appeared  in \cite{cklmns12}.   We refer to \cite{cc12} for some other results about scalable frames.

\subsection{Probabilistic frames}\label{subsec1.3}  While frames are intrinsically defined through their spanning properties, in real euclidean spaces, they   can also be viewed as  distributions of point masses. In this context, the notion of \emph{probabilistic frames}\index{probabilistic frame} was introduced as a class of probability measures with finite second moment and whose support spans the entire space \cite{ehloko10, ehler10, ehloko11}. Probabilistic frames   are special cases of continuous frames\index{continuous frames} as introduced by S.~T.~Ali, J.-P.~Antoine, and P.-P.~Gazeau \cite{aag}, see also \cite{fora}.

In Section~\ref{sec3}, we consider frames from this probabilistic point of view.  To begin we note  that probabilistic frames are extensions of the notion of frames previously defined. Indeed,  consider  a frame, $\Phi=\{\varphi_k\}_{k=1}^M $ for $\R^N$ and define the probability measure $$\mu_\Phi = \tfrac{1}{M}\sum_{k=1}^M \delta_{\varphi_{k}}$$ where $\delta_x$ is the Dirac mass at $x \in \R^N$. It is easily seen that the second moment\index{second moment} of $\mu_\Phi$ is finite, i.e., $$\int_{\R^N}\|x\|^2d\mu_{\Phi}(x)= \tfrac{1}{M}\sum_{k=1}^M \|\varphi_k\|^2<\infty,$$ and the span of the support of $\mu_\Phi$ is $\R^N$. Thus, each (finite) frame can be associated to a probabilistic frame. We shall present other examples of probabilistic frames associated to $\Phi$ in Section~\ref{sec3}.  
By analogy to the theory of finite frame,  we shall say that a probability measure on $\R^N$ with finite second moment is a \emph{probabilistic frame} if the linear span of its support is $\R^N$.  

It is known that  the space $\mathcal{P}_2(\R^N)$ of probability measures with finite second moments  can be equipped with the ($2$-)Wasserstein metric\index{Wasserstein metric}. In this setting, many questions  in frame theory can  be seen as analysis problems on a subset of the Wasserstein metric space  $\mathcal{P}_2(\R^N)$. In Section~\ref{sec3} we introduce this metric space and derive  some immediate properties of frames in this setting. Moreover, this probabilistic setting allows one to use non-algebraic tools to investigate questions from frame theory. 
For instance, tools like convolution of measures has been used in \cite{ehloko11} to build new (probabilistic) frames from old ones.

One of the main advantages in analyzing frames in the context of the Wasserstein metric lies in the powerful tools available to solve some optimization problems involving frames using the  framework of optimal transport. While we will not delve into any details in this chapter, we point out that optimization of functionals like the frame potential can be studied in this setting. For example, C.~Wickman recently showed that a potential function that generalizes Benedetto and Fickus's frame potential  can be minimized in the Wasserstein space using  some optimal transport notions \cite{cw14}. 
In the last part of the lecture, we shall focus on a family of potentials that generalize the frame potential and present a survey of recent results involving their minimization. In particular, this family includes  the \emph{coherence} of a set of vectors, which is important quantity in compressed sensing \cite{strhea, wal03, welch74}, as well as a functional whose minimizers are, in some cases, solutions to the Zauner's conjecture \cite{renes, zau99}. 

Though we shall not elaborate on these here, it is worth mentioning that probabilistic frames are related to many other areas including: (a) the covariance of matrices  multivariate random vectors \cite{elrv12, mr95, mr97, rv12, rv01}; (b)  directional statistics where there are used to test whether certain data are uniformly distributed; see, \cite{mejg11, Mardia:2008aa}; (c)  isotropic measures \cite{gimi00, Milman:1987aa}, which, as we shall show, are related to the class of \emph{tight probabilistic frames}. We refer to \cite{ehloko11} for an overview of other relationships between probabilistic frames and other areas.

The rest of this chapter  is organized as follows. In Section~\ref{sec2} we give an overview of the recent developments on scalable frames or preconditioning of frames. In Section~\ref{sec3} we deal with probabilistic frames and outline the recently introduced Wasserstein metric tools to deal with frame theory.

\section{Preconditioning techniques in frame theory}\label{sec2}

Scalable frames\index{scalable frame} are frames $\Phi=\{\varphi_k\}_{k=1}^M$ for $ \R^N$ for which there exist nonnegative scalars $\{c_k\}_{k=1}^M \subset [0, \infty)$ such that $$\tilde{\Phi}=\{c_k \varphi_k\}_{k=1}^M$$ is a tight frame. 
Scalable frames were first introduced and characterized in \cite{kopt12}. Both infinite and finite dimensional settings were considered. In this section, we only focus on the latter giving an overview of recent methods developed  to understand  scalable frames.  The results that we shall describe give an exact characterization of the set of scalable frames from various perspectives. However,  the important and very practical question of developing algorithms to find the weights $\{c_k\}$ that make the frame tight will not be considered here, and we refer to \cite{ccko15} for a sample of results on this topic. Similarly, when a frame fails to be scalable, one could seek to relax the tightness condition and seek an ``almost scalable frame''. These considerations are sources of ongoing research and will not be taken upon here.  Finally, it is worth pointing out a very interesting application of the theory of scalable frames to wavelets constructed from the Laplacian Pyramid scheme \cite{yhko15}. 

The rest of this section is organized as follows. In Section~\ref{subsec2.1} we define scalable frames and derive some of their elementary properties. We then outline a characterization of the set of scalable frames in terms of certain convex polytopes in Section~\ref{subsec2.2}. This characterization is preceded by motivating examples of scalable frames in dimension $2$.  In Sections~\ref{subsec2.3} and~\ref{subsec2.4} we give two other equivalent characterizations of scalable frames. The first of these characterizations has geometric  interpretation, while the second one is based on Fritz John's ellipsoid theorem.

\subsection{Scalable frames: Definition and properties }\label{subsec2.1} The following definitions of scalable frames first appeared in \cite{kopt12, kop13c}: 

\begin{definition}\label{def1}
Let $ M, N$ be integers  such that $N  \leq  M$. A frame $\Phi=\{\varphi_k\}_{k=1}^{M}$ in $\R^N$ is called  \emph{scalable}\index{scalable frame}, \emph{respectively, strictly scalable,}\index{scalable frame!strictly scalable frame} if there exist nonnegative, respectively, positive, scalars $\{c_k\}_{k=1}^M$ such that $\{c_k\varphi_k\}_{k=1}^M$ is a tight frame for $\R^N$. The set of  scalable frames, respectively,  strictly scalable frames,  is denoted by $\sca(M,N)$, respectively,  $\sca_{+}(M,N)$. 

Moreover,  given an integer $m$ with $N\leq m \leq M$, $\Phi=\{\varphi_k\}_{k=1}^{M}$ 
is said to be \emph{$m$-scalable}\index{scalable frame!$m$-scalable}, respectively,  \emph{strictly $m-$scalable}\index{scalable frame!strictly $m$-scalable}, if there exist a subset $\Phi_I=\{\varphi_k\}_{k\in I}$ with $I\subseteq \{1, 2, \hdots, M\}$, $\#I=m$,  such that $\Phi_I=\{\varphi_k\}_{k\in I}$ is scalable, respectively, strictly scalable. 

We denote the set of $m$-scalable frames, respectively, strictly $m$-scalable frames  in $\fr(M,N)$ by   $\sca(M,N,m)$, respectively, $\sca_{+}(M,N,m).$ 
\end{definition}

When the integer $m$ is fixed in a given context, we will simply refer to an $m-$scalable frame as a \emph{scalable frame.}  The role of the parameter $m$ is especially relevant when dealing with frames of very large redundancy, i.e., when $M/N \gg 1$. In such a case, choosing a ``reasonable'' $m$ such that the frame is $m-$scalable could potentially lead to sparse representations for signals in $\R^N$. In addition, the problems of finding the weights that make a frame scalable as well as determining the  smallest $m$ such that a given frame is $m-$ scalable have been considered in \cite{ccnst14, cc12}. We shall give more details about this question in Section~\ref{subsec2.2}.

We now point out some special and trivial examples of scalable frames. When $M=N$, a frame $\Phi$ is scalable if and only if $\Phi$ is an orthogonal set. Moreover, when $M\geq N$, if $\Phi$ contains an orthogonal basis, then it is clearly $N-$scalable. Thus, given $M\geq N$, the set $\sca(M,N,N)$ consists exactly of frames that contains an orthogonal basis for $\R^N$. 

So from now on we shall assume without loss of generality that $M\geq N+1$, that $\Phi$ contains no orthogonal basis, and that $\varphi_k \neq \pm \varphi_\ell$ for $\ell\neq k$.

Given a frame $\Phi \subset \R^N$, assume that $\Phi=\Phi_1 \cup \Phi_2$ where 
$$\Phi_1=\{\varphi_{k}^{(1)} \in \Phi: \varphi_{k}^{(1)}(N)\geq 0\}$$  and $$\Phi_2=\{\varphi_{k}^{(2)} \in \Phi: \varphi_{k}^{(2)}(N)< 0\}.$$ In other words, $\Phi_1$ consists of all frame vectors from $\Phi$ whose $N^{th}$ coordinates are nonnegative. Then the frame $\Phi'=\Phi_1 \cup (-\Phi_2)=\{\varphi_{k}^{(1)}\}\cup\{-\varphi_{k}^{(2)}\}$ has the same frame operator as $\Phi$. In particular, $\Phi$  is a tight frame if and only if $\Phi'$ is a tight frame.   In addition, $\Phi$ is scalable if and only if $\Phi'$ is scalable with exactly the same set of weights. Note that the frame vectors in $\Phi'$ are all in the upper-half space. Thus, when convenient we shall assume without loss of generality that all the frame vectors are in the upper-half space,  that is $\Phi \subset \R^{N-1}\times \R_{+}$ where $\R_{+}=[0, \infty)$.

We note  that a frame $\Phi=\{\varphi_k\}_{k=1}^M\subset \R^N$  with $\varphi_k \neq 0$ for each $k=1, \hdots, M$ is scalable if and only if $\Phi'=\{\tfrac{\varphi_k}{\|\varphi_k\|}\}_{k=1}^M$ is scalable. Consequently, we might assume in the sequel that we work  with frames consisting of unit norm vectors.

We now collect a number of elementary properties of the set of scalable frames in $\R^N$.  We refer to  \cite{kopt12, kop13c} for details. 

\begin{proposition}\label{proposition:recaps} Let $M\geq N$, and $m\geq 1$ be integers. 
\begin{enumerate}
\item[(i)] If $\Phi\in\fr$ is  $m$-scalable then $m\geq N$. 
\item[(ii)] For any integers $m, m'$ such that  $N\leq m\leq m'\leq M$ we have that $$\sca(M,N, m)\subset \sca(M,N, m'),$$ and  $$
\sca(M,N) = \sca(M,N,M) = \bigcup_{m=N}^{M}\sca(M,N, m).
$$
\item[(iii)] $\Phi\in\sca(M,N)$  if and only if $T(\Phi)\in\sca(M,N)$ for one (and hence for all) orthogonal
transformation(s) $T$ on $\R^N$.
\item[(iv)] Let $\Phi=\{\varphi_k\}_{k=1}^{N+1} \in\fr(N+1, N) \setminus \{0\}$ with $\varphi_k\neq \pm \varphi_\ell$ for $k\neq \ell$. If $\Phi \in \sca_{+}(N+1, N,N)$, then $\Phi \notin \sca_{+}(N+1,N)$. 
\end{enumerate}
\end{proposition}
 
\begin{remark}
We point out that part (iii) of Proposition~\ref{proposition:recaps} is equivalent to saying that $\Phi$ is not scalable if one can find an orthogonal transformation $T$ on $\R^N$ such that $T(\Phi)$ is not scalable. 
\end{remark}

Besides these elementary properties, a study of the topological properties of the set of scalable frames was considered in \cite{kopt12, kop13c}. In particular, 
\begin{proposition}\label{proposition:topo} Let $M\geq N\geq 2$.
\begin{enumerate}
\item[(i)] $\sca(M,N)$ is closed in $\fr(M,N)$. Furthermore, for each $N\leq m \leq M$, $\sca(M,N,m)$ is closed in $\fr(M,N)$.
\item[(ii)] If $M< N(N+1)/2$, then the interior of $\sca(M,N)$ is empty.
\end{enumerate}
\end{proposition}

Part (i) of proposition~\ref{proposition:topo} was proved in \cite[Corrollary 3.3]{kopt12} and \cite[Proposition 4.1]{kop13c}, while part (ii) first appeared in \cite[Theorem 4.2]{kop13c}.

\subsection{Convex polytopes associated to  scalable frames }\label{subsec2.2}
We now proceed to write an explicit formulation of the scalability problem.  From this formulation a convex geometric characterization of $\sca(M,N)$ will follow. To start, we  recall that $\Phi$ denote the synthesis operator associated to the frame $\Phi=\{\varphi_k\}_{k=1}^{M}$.  $\Phi$ is ($m$-) scalable if and only if there are  positive numbers $\{x_k\}_{k\in I} $  with $\#I=m \geq N$   such that $\widetilde{\Phi}=\Phi X$  satisfies
\begin{equation}\label{eqscfrm}
\widetilde{\Phi} \widetilde{\Phi}^{T}=\Phi X^{2} \Phi^{T}=\tilde{A}I_{N}=\tfrac{\sum_{k\in I}x_{k}^{2}\|\varphi_{k}\|^{2}}{N}I_{N}
\end{equation} where $X$ is the diagonal matrix with the weights $x_k$ on its main diagonal if $k\in I$ and $0$  for $k\notin I$, and $I_{N}$ is the $N\times N$ identity matrix. Moreover, $$\|X\|_{0}=\#\{k: x_k>0\}=m \geq N.$$ By rescaling the diagonal matrix $X$, we can assume that $\tilde{A}=1$.  Thus,~\eqref{eqscfrm} is equivalent to solving 
\begin{equation}\label{eqscfrm-reduced}
\Phi Y\Phi^{T}=I_{N}
\end{equation} for $Y=\tfrac{1}{\tilde{A}} X^2.$

To gain some intuition let us consider the two dimensional case with $M\geq 3$. In particular, let us describe when  $\Phi=\{\varphi_k\}_{k=1}^M \subset S^{1}$ is  a   scalable frame.   Without loss of generality, we may assume that $\Phi=\{\varphi_{k}\}_{k=1}^{M} \subset \R \times \R_{+}$,  $\|\varphi_k\|=1$, and $\varphi_\ell \neq \pm \varphi_k$ for $\ell \neq k$.  Thus $$\varphi_{k}=\begin{pmatrix} \cos \theta_k\\ \sin \theta_k\end{pmatrix} \in S^{1}$$with $$0=\theta_1< \theta_2 < \theta_3 < \hdots < \theta_M <\pi.$$  Let $Y=(y_k)_{k=1}^{M}\subset [0, \infty)$, then~\eqref{eqscfrm-reduced} becomes

\begin{equation}\label{2-dM}
\begin{pmatrix} \sum_{k=1}^M y_k\cos^2 \theta_k & \sum_{k=1}^M y_k\sin\theta_k \cos \theta_k\\ 
\sum_{k=1}^M y_k\sin\theta_k \cos \theta_k & \sum_{k=1}^M y_k\sin^2\theta_k \end{pmatrix}=\begin{pmatrix}1 & 0 \\ 0 & 1\end{pmatrix}.
\end{equation} This is equivalent to 

$$ \left\{ \begin{array} {r@{\quad = \quad}l}
\sum_{k=1}^M y_k\cos^2 \theta_k &1\\
\sum_{k=1}^M y_k\sin^2\theta_k&1\\
\sum_{k=1}^M y_k\sin\theta_k \cos \theta_k &0,
\end{array}\right.$$ and using some row operations we arrive at

$$ \left\{ \begin{array} {r@{\quad = \quad}l}
\sum_{k=1}^M y_k\sin^2 \theta_k &1\\
\sum_{k=1}^M y_k\cos2\theta_k &0\\
\sum_{k=1}^M y_k\sin2\theta_k  &0.
\end{array}\right.$$ For $\Phi$ to be scalable we must find a nonnegative vector $Y=(y_k)_{k=1}^M$ in the kernel of the matrix whose $k^{th}$ column is $\begin{pmatrix}\cos2\theta_k\\ \sin2\theta_k\end{pmatrix}.$ Notice that the first equation is just a normalization condition. 

 We now describe the the subset of the kernel of this matrix that consists of non-trivial nonnegative vectors. 
Observe that the matrix  can  be reduced to 

\begin{equation}\label{3redumatbis}
\begin{pmatrix} 1 & \cos 2\theta_2 & \hdots & \cos 2\theta_M \\ 0 & \sin 2\theta_2 & \hdots & \sin 2\theta_M\end{pmatrix}.
\end{equation}

\begin{example}\label{example2-3}
We first consider the case $M=3$. In this case, we have $0=\theta_1< \theta_2 < \theta_3 < \pi,$
and the~\eqref{3redumatbis} becomes

\begin{equation}\label{3redumat3}
\begin{pmatrix} 1 & \cos 2\theta_2 & \cos 2\theta_3 \\ 0 & \sin 2\theta_2 &  \sin 2\theta_3\end{pmatrix}.
\end{equation}

  If  there exists an index $k_0 \in \{2, 3\}$ with $\sin 2\theta_{k_{0}}=0$, then $\theta_{k_{0}}=\pi/2$ and the corresponding frame contains an ONB and, hence is scalable.  
\begin{enumerate}
\item[(i)] Moreover, if  $k_0=2$, then $0=\theta_1<\theta_2=\pi/2 <\theta_3<\pi $.  In this case, the fame is $2-$ scalable but not $3-$ scalable, i.e.,  the frame is in $\sca_{+}(3, 2, 2) \setminus \sca(3, 2, 3)$. This is illustrated by  Figure~\ref{fig-1a}. 
\item[(ii)]  If  $k_0=3$, then $0=\theta_1<\theta_2<\theta_3=\pi/2$.  By symmetry (with respect to the $y$ axis) we conclude again that the fame is $2-$ scalable but not $3-$ scalable.
\end{enumerate}

\begin{figure}[h]
\begin{center}
\includegraphics[scale=0.33]{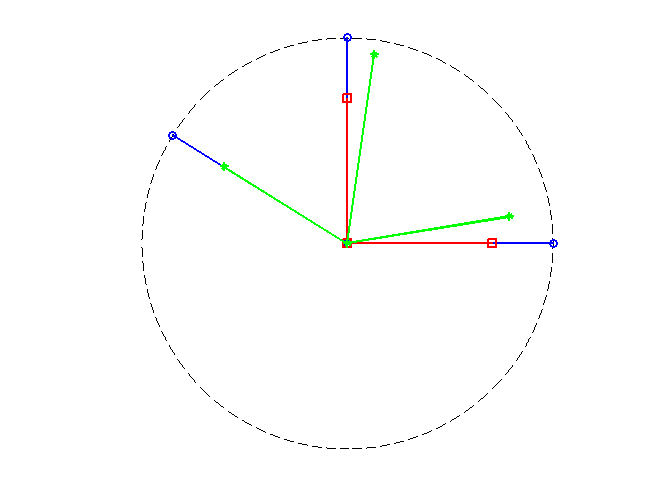} 
\end{center}
\caption{A scalable frame (contains an orthonormal basis) with $3$ vectors in $\R^2$. The original frame vectors are in blue, the frame vectors obtained by scaling  are in red, and for comparison the associated canonical tight frame vectors are  in green. }
\label{fig-1a}
\end{figure}

Assume now that  $\theta_k\neq \pi/2$ for $k=2, 3$. If $\theta_3< \pi/2$, then the frame cannot be scalable. Indeed, $u=(z_1, z_2, z_3)$ belongs to the kernel of~\eqref{3redumat3} if and only if 

\begin{equation}\label{solution2-3}
\left\{ \begin{array} {r@{\quad = \quad}l}
z_1&\tfrac{\sin 2(\theta_3-\theta_2)}{\sin 2\theta_2}z_3,\\
z_2&-\tfrac{\sin 2\theta_3}{\sin 2\theta_2}z_3,
\end{array}\right.
\end{equation} where $z_3\in \R$. 
The choice of the angles implies that $z_2 z_3 \leq 0$ and $z_1z_3\leq 0$ with equality if and only if  $z_3=0$.  This is illustrated by  Figure~\ref{fig-1b}.
Similarly, if $0=\theta_1< \pi/2 < \theta_2< \theta_3<\pi$, then the frame cannot be scalable. 

\begin{figure}[h]
\begin{center}
\includegraphics[scale=0.33]{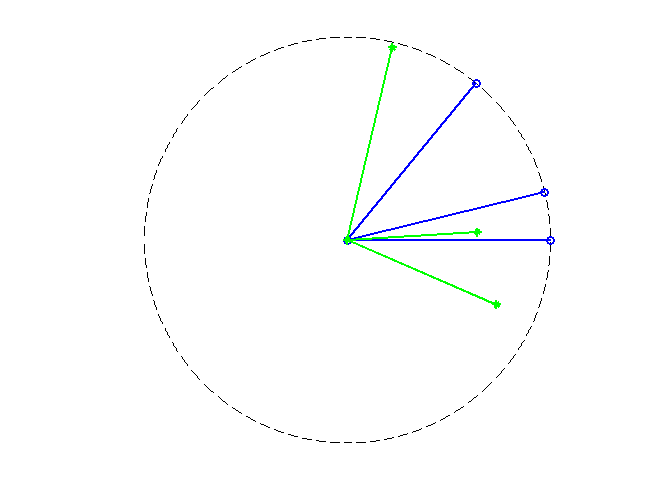}
\end{center}
\caption{A non scalable frame with $3$ vectors in $\R^2$. The original frame vectors are in blue,  for comparison the associated canonical tight frame vectors are  in green.}
\label{fig-1b}
\end{figure}

 On the other hand  if $0=\theta_1< \theta_2< \pi/2 < \theta_3<\pi$, then it follows  from~\eqref{solution2-3} $z_1>0$ and $z_2>0$  for all $z_3>0$  if and only if $\theta_3-\theta_2<\pi/2$. Consequently,
  when $0=\theta_1< \theta_2< \pi/2 < \theta_3<\pi$ the frame $\Phi \in \sca_+(3,2,3)$ if and only if $0< \theta_3-\theta_2<\pi/2$.   This is illustrated by Figure~\ref{fig-1c}.

\begin{figure}[h]
\begin{center}
\includegraphics[scale=0.33]{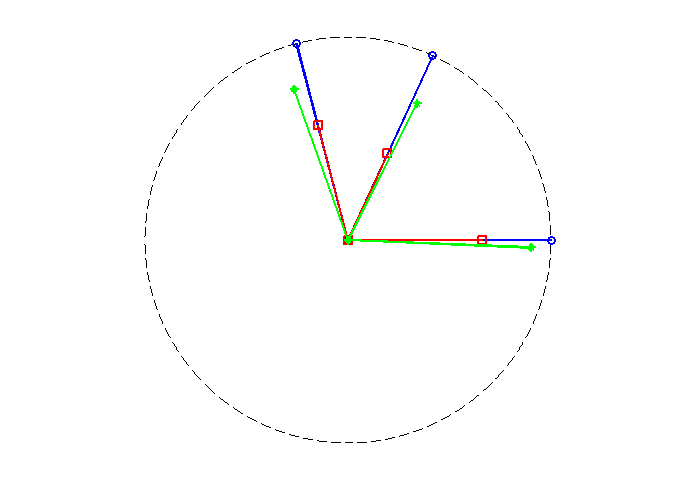}
\end{center}
\caption{A scalable frame  with $3$ vectors in $\R^2$. The original frame vectors are in blue, the frame vectors obtained by scaling  are in red, and for comparison the associated canonical tight frame vectors are  in green.}
\label{fig-1c}
\end{figure}

\end{example}

\begin{example}\label{example2-4}
Assume now that $M=4$. Then we are lead to seek nonnegative non-trivial vectors in the null space of
\begin{equation*}
\begin{pmatrix} 1 & \cos 2\theta_2 & \cos 2\theta_3 & \cos 2\theta_4 \\ 0 & \sin 2\theta_2 &  \sin 2\theta_3 &  \sin 2\theta_4\end{pmatrix}.
\end{equation*}

If there exists an index $k_0 \in \{2, 3, 4\}$ with $\sin 2\theta_{k_{0}}=0$, then $\theta_{k_{0}}=\pi/2$ and the corresponding frame contains an ONB. Consequently, the frame is scalable.  In particular, 
\begin{enumerate}
\item When $k_0=2$,  the null space of the matrix is described by 
\begin{equation*}
\left\{ \begin{array} {r@{\quad = \quad}l}
z_1&z_2 + \tfrac{\sin 2(\theta_4-\theta_3)}{\sin 2\theta_3}z_4,\\
z_3&-\tfrac{\sin 2\theta_4}{\sin 2\theta_3}z_4,
\end{array}\right.
\end{equation*} where $z_2, z_4\in \R$.  Note that  $z_3\leq 0,$ with equality only when $z_4=0$, in which case $z_3=0$ and the frame will be $2-$scalable, but not $m-$scalable for $m=3, 4$. This is illustrated by the left figure in   Figure~\ref{fig-4a}.

\begin{figure}[h]
\begin{center}
\includegraphics[scale=0.3]{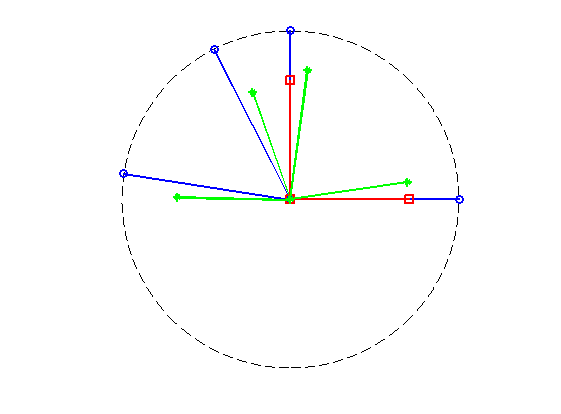}
\end{center}
\caption{ A scalable frame (contains an orthonormal basis) with with $4$ vectors in $\R^2$. The original frame vectors are in blue, the frame vectors obtained by scaling  are in red, and for comparison the associated canonical tight frame vectors are  in green.}
\label{fig-4a}
\end{figure}

\item If instead, $k_0=3$,  then a similar argument shows that 
\begin{equation*}
\left\{ \begin{array} {r@{\quad = \quad}l}
z_1&z_3 + \tfrac{\sin 2(\theta_4-\theta_2)}{\sin 2\theta_2}z_4,\\
z_2&-\tfrac{\sin 2\theta_4}{\sin 2\theta_2}z_4,
\end{array}\right.
\end{equation*} where $z_3, z_4\in \R$. Any choice of $z_4>0$ will result in $z_2>0$. If we choose $\theta_4-\theta_2<\pi/2$, then $z_3\geq 0$ will lead to a $3-$ scalable frame or a $4-$scalable frame. If instead, $\theta_4-\theta_2\geq \pi/2$ we can always choose $z_3>0$ and large enough to guarantee that $z_1>0$, hence $\Phi$ will be $4-$scalable.  
\item When $k_0=4$, then $\Phi \in \sca_{+}(4, 2,2) \setminus \sca(4, 2, m)$ for $m =3, 4$.   
\item When $\sin 2\theta_k\neq 0$ for $k \in \{2, 3, 4\}$ then
\begin{equation*}
\left\{ \begin{array} {r@{\quad = \quad}l}
z_1&\tfrac{\sin 2(\theta_3-\theta_2)}{\sin 2\theta_2}z_3 + \tfrac{\sin 2(\theta_4-\theta_2)}{\sin 2\theta_2}z_4,\\
z_2&-\tfrac{\sin 2\theta_3}{\sin 2\theta_2}z_3-\tfrac{\sin 2\theta_4}{\sin 2\theta_2}z_4,
\end{array}\right.
\end{equation*} where $z_3, z_4\in \R$. A choice of $z_3, z_4\geq 0$ will lead to a scalable frame if at least $z_1\geq 0$ or $z_2\geq 0$.  For example, Figure~\ref{fig-4b} shows a scalable frame. 

\begin{figure}[h]
\begin{center}
\includegraphics[scale=0.3]{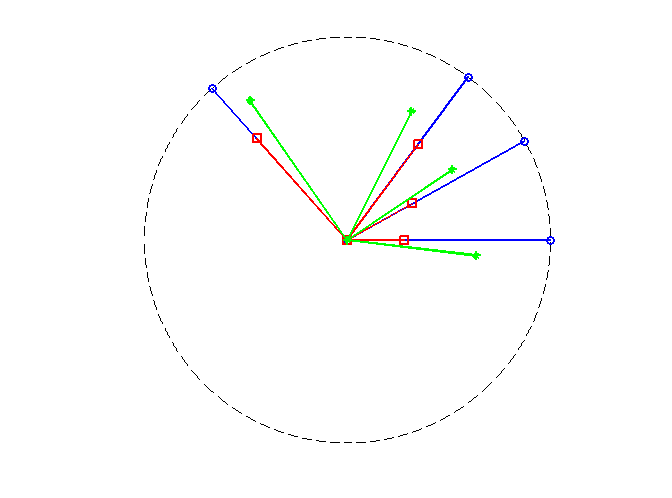}
\end{center}
\caption{A scalable frame with with $4$ vectors in $\R^2$. The original frame vectors are in blue, the frame vectors obtained by scaling  are in red, and for comparison the associated canonical tight frame vectors are  in green.}
\label{fig-4b}
\end{figure}

But in this case we could also get non scalable frame, see  Figure~\ref{fig-4c} . The implication here is that the scalability of the frames depends on  the relative position of the frame vectors (hence) the angles $\theta_k$. This will be made rigorous in  Section~\ref{subsec2.3}.
\end{enumerate}

\begin{figure}[h]
\begin{center}
\includegraphics[scale=0.3]{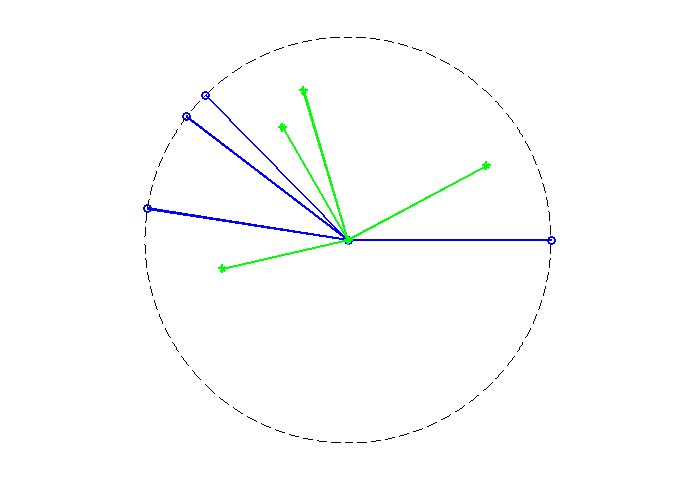}
\end{center}
\caption{A non scalable frame  with $4$ vectors in $\R^2$. The original frame vectors are in blue,  for comparison the associated canonical tight frame vectors are  in green.}
\label{fig-4c}
\end{figure}

\end{example}

More generally in this two dimensional case, we can continue this analysis  of the transformation given by the matrix~\eqref{3redumatbis} to characterize when $\Phi=\{\varphi_k\}_{k=1}^M$ is scalable. From the figures shown in Examples~\ref{example2-3} and~\ref{example2-4},  it is clear that some geometric considerations are involved.
Before elaborating more on these geometric considerations in Section~\ref{subsec2.3}, we now consider the general case $M\geq N\geq 2$.  In particular, we follow the algebraic approach given in the two-dimensional case, by  writing out the equations in~\eqref{eqscfrm-reduced} and collecting all the diagonal terms one the one hand, and the non-diagonal terms on the other, we see that for a frame to be
$m$-scalable it is necessary and sufficient that there exists $u = (c_{1}^{2},c_{2}^{2},\ldots,c_{M}^{2})^{T}$
with $\|u\|_0:=\# \{u_k: u_k>0\} \leq m$ which is a solution of the following linear system of $\tfrac{N(N+1)}{2}$ equations in the $M$ unknowns $(y_j)_{j=1}^M$:
\begin{equation}\label{nxm}
\left\{ \begin{array} {r@{\quad = \quad}l}
\sum\limits_{j=1}^M\varphi_j(k)^2y_j & 1 \quad \text{for }k=1,\ldots,N,\\
\sum\limits_{j=1}^M\varphi_j(\ell)\varphi_j(k)y_j & 0\quad  \text{for }k>\ell =1,\ldots,N.
\end{array}\right.
\end{equation}

We can further reduce this linear system in the following manner. We keep all the equations with homogeneous right-hand sides, i.e., those coming from the non diagonal terms of~\eqref{eqscfrm-reduced}. There are $N(N-1)/2$ such equations. The remaining $N$ equations come from the diagonal terms of~\eqref{eqscfrm-reduced}, and their right hand-sides are all  $1$. We can use row  operations to reduce these to a new set of $N$ linear equations the first of which will be $$\sum_{j=1}^M\varphi_j(1)^2y_j=1.$$ For $k=2, \hdots, N$, the $k^{th}$ equation is obtained by subtracting row $1$ from row $k$ leading to $$\sum_{j=1}^{M}(\varphi_{j}(k)^2-\varphi_{j}(1)^2)y_j=0.$$ The condition $$\sum_{j=1}^M\varphi_j(1)^2y_j=1$$ is a normalization condition, indicating that if $\Phi$ can be scaled with $y=(y_j)_{j=1}^M\subset [0, \infty)$, then it can also be scaled by $\lambda y $ for any $\lambda>0$.  Thus, ignoring this condition and collecting all the remaining equations, we see that $\Phi$ is $m$-scalable if and only
if there exists a nonnegative vector $u\in\R^M$ with $\|u\|_0\le m$  such that $$F(\Phi)u=0,$$ where 
the $(N-1)(N+2)/2\times M$ matrix $F(\Phi)$ is given by
\begin{equation*}
F(\Phi) = \begin{pmatrix} F(\varphi_{1}) & F(\varphi_{2})  & \hdots & F(\varphi_{M})  \end{pmatrix},
\end{equation*}
where $F : \R^N\to\R^{d}$, $d := (N-1)(N+2)/2$, is defined by
\begin{equation}\label{fctF}
F(x)=
\begin{pmatrix}
F_{0}(x)\\ F_{1}(x)\\ \vdots\\ F_{N-1}(x)
\end{pmatrix},
\qquad
F_{0}(x) =
\begin{pmatrix} x_{1}^{2}-x_{2}^{2}\\ x_{1}^{2}-x_{3}^{2} \\ \vdots \\ x_{1}^{2}-x_{N}^{2}\end{pmatrix},
\qquad
F_{k}(x) =
\begin{pmatrix}x_{k}x_{k+1}\\ x_{k}x_{k+2}\\ \vdots \\ x_{k}x_{N}\end{pmatrix},
\end{equation}
and $F_{0}(x)\in \R^{N-1}$, $F_{k}(x) \in \R^{N-k}$, $k=1, 2, \hdots, N-1$.

To summarize, we arrive at the following result  that was proved in \cite[Proposition 3.7]{kop13c}

\begin{proposition}\cite[Proposition 3.7]{kop13c}\label{pro:nullspaceFPHI}
A frame $\Phi$ for $\R^N$ is  $m$-scalable, respectively,  strictly $m$-scalable,  if and only if there exists a nonnegative
$u\in\ker F(\Phi)\setminus\{0\}$ with $\|u\|_0\le m$, respectively,  $\|u\|_0 = m$, where $\|u\|_0=\#\{k: u_k>0\}.$
\end{proposition}

In the two dimensional case the map $F$ reduces to 
\begin{equation*}
F\begin{pmatrix} x\\ y\end{pmatrix}=\begin{pmatrix} x^2 -y^2\\ xy\end{pmatrix}
.\end{equation*}
However, in all  the previous examples  we considered instead the more geometric map : 
\begin{equation*}
\widetilde{F}\begin{pmatrix} x\\ y\end{pmatrix}=\begin{pmatrix} x^2 -y^2\\ 2xy\end{pmatrix}
.\end{equation*}
 It is readily seen that $F(\Phi)$ and $\widetilde{F}(\Phi)$ have exactly the same kernel. In fact the map $\widetilde{F}$ carries the following geometric interpretation.  Let $L_{\theta}$ be a line through the origin in $\R^2$ which makes an angle $\theta$ with the positive $x-$axis with $\theta \in [0, \pi]$. Then the image of $L_{\theta}$ by $\widetilde{F}$ is the line $L_{2\theta}$ that makes an angle $2\theta$ with the positive $x-$ axis. That is, $\widetilde{F}$ just rotates counterclockwise the line $L_{\theta}$ by an angle equal to $\theta$.

In the two dimensional case we exploited the geometric meaning of the map $F$ or $\widetilde{F}$ to describe the subset of nonnegative vectors of the nullspace of  $\widetilde{F}(\Phi)$.  More generally, to find nonnegative vectors in the nullspace of the matrix $F(\Phi)$ we can appeal to one of the formulation of Farkas lemma\index{Farkas lemma}: 

\begin{lemma}\cite[Lemma  1.2.5]{matou02}\label{l:farkas}
For every real $N\times M$-matrix $A$ exactly one of the following cases occurs:
\begin{itemize}
\item[(i)] The system of linear equations $Ax=0$ has a nontrivial nonnegative solution $x\in \R^{M}$, i.e., all components of $x$ are nonnegative and at least one of them is strictly positive.
\item[(ii)] There exists $y\in\R^N$ such that $y^{T}A$ is a vector with all entries strictly positive.
\end{itemize}
\end{lemma}

Applying this in the two dimensional case, we see that for the frame to be scalable, the second alternative in Farkas's lemma should not hold. That is there must exist no vector in $\R^2$ that lies on ``one side'' of all the vectors $F(\varphi_k)$
for $k=1, 2, \hdots, M$. We illustrate this by the following figures:

\begin{figure}[h]
\begin{center}
\includegraphics[scale=0.65]{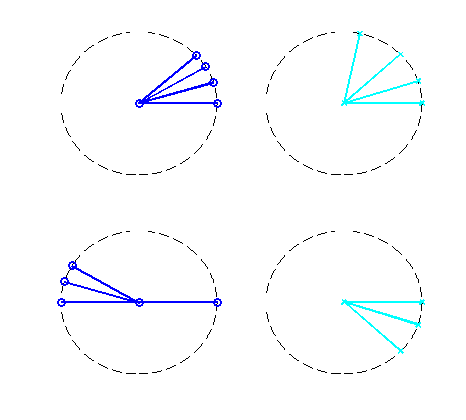}
\end{center}
\caption{Frames with $4$ vectors in $\R^2$ (in blue, top and bottom left) and their images by the map $F$ (in green, top and bottom right). Both of these examples result in non scalable frames.}
\label{fig-3}
\end{figure}

\begin{figure}[h]
\begin{center}
\includegraphics[scale=0.65]{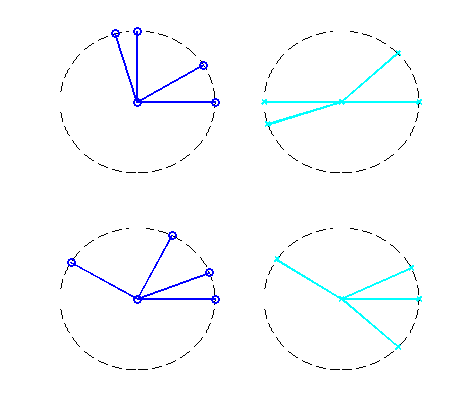}
\end{center}
\caption{Frames with $4$ vectors in $\R^2$ (in blue, top and bottom left) and their images by the map $F$ (in green, top and bottom right). Both of these examples result in scalable frames.}
\label{fig-4}
\end{figure}

We make the following observation about the first alternative of Lemma~\ref{l:farkas}. If $\{A_k\}_{k=1}^M\subset \R^N$ represents  the column vectors of $A$, then there exists of a vector $0\neq x=(x_k)_{k=1}^M$ with $x_k\geq 0$ such that $Ax=0$ is equivalent to saying that $$\sum_{k=1}^Mx_kA_k=0.$$ Without loss of generality we may assume that $\sum_{k=1}^Mx_k=1$, in which case the condition is equivalent to $0$ being a convex combination of the column vectors of $A$. Thus, having a nontrivial nonnegative vector in the null space of $A$ is a statement about the convex hull of the columns of $A$.

Motivated by the geometric intuition we gained from the two-dimensional setting and to  effectively use Farkas's lemma, we introduce a few notions from convex geometry, especially the theory of convex polytopes\index{convex polytope}, and we refer to \cite{stwi70, web94} for more details on these concepts. For a finite set  $X =\{x_i\}_{k=1}^{M} \subset \R^N$, the \emph{polytope generated by $X$} is the 
convex hull of $X$, which is  a  compact convex subset of $\R^N$. In particular, we denote this set by  $P_{X}$ (or $\co(X)$), and we have $$
P_{X}=\co(X):= \left\{ \sum_{k=1}^{M} \alpha_{k} x_k:  \alpha_k \geq 0, \, \sum_{k=1}^{M}\alpha_k=1\right\}.
$$
The \emph{affine hull generated by $X$} is defined by
$$
\aff(X) := \left\{ \sum_{k=1}^{M} \alpha_{k} x_k:   \, \sum_{k=1}^{M}\alpha_k=1\right\}.
$$
We have $\co(X) \subset\aff(X)$. 
The \emph{relative interior}\index{relative interior} of the polytope $\co(X)$ denoted by $ri\co(X)$, is the interior of $\co(X)$ in the
topology induced by $\aff(X)$. We have that $ri\co(X)\neq\emptyset$ as long as $\#X\ge 2$, and 

$$
ri\co(X) = \left\{\sum_{k=1}^{M}\alpha_k x_k: \alpha_k>0, \sum_{k=1}^{M}\alpha_k = 1\right\}
.$$

The \emph{polyhedral cone generated by $X$} is the closed convex cone\index{cone} $C(X)$ defined by
$$
C(X) = \left\{\sum_{i=k}^{M} \alpha_{k} x_k : \alpha_k\geq 0\right\}.
$$
 The \emph{polar cone} of $C(X)$ is the  closed convex cone $C^{\circ}(X)$ defined by
$$
C^{\circ} (X):= \{x\in\R^N : \ip{x}{y}\leq 0\,\,{\textrm for \, all}\,\,y\in C(X)\}.
$$
The cone $C(X)$ is said to be \emph{pointed} if $C(X)\cap (-C(X)) = \{0\},$ and \emph{blunt} if the linear space generated
by $C(X)$ is $\R^N$, i.e.,  $\linspan C(X) = \R^N$.

Using Proposition~\ref{pro:nullspaceFPHI}, we see that $\Phi$ is ($m-$)scalable if there exists $\{\lambda_k\}_{k\in I}\subset [0, \infty)$ , $\# I =m$ such that $$\sum_{k\in I}\lambda_k F(\varphi_k)=0.$$ This is equivalent to saying that $0$ belongs to the polyhedral cone generated by $F(\Phi_I)=\{F(\varphi_k)\}_{k\in I}.$ Without loss of generality we can assume that $\sum_{k\in I}\lambda_k=1$ which implies that $0$ belongs to the polytope generated by $F(\Phi_I)=\{F(\varphi_k)\}_{k\in I}.$ Putting these observations together with Lemma~\ref{l:farkas}  the following results were established  in \cite[theorem 3.9]{kop13c}. In the sequel, we shall denote by $[K]$ the set of integers  $\{1, 2, \hdots, K\}$ where  $K \in \N$. 

\begin{theorem}\cite[Theorem 3.9]{kop13c}\label{poly}
Let $M\geq N\geq 2$, and let $m$ be such that $N\leq m\leq M$. Assume that $\Phi=\{\varphi_k\}_{k=1}^M \in \fr^{*}(M,N)$ is such that  $\varphi_k \neq \varphi_\ell$ when $k\neq \ell$. Then the following statements are equivalent:
\begin{enumerate}
\item[(i)] $\Phi$ is $m-$scalable, respectively, strictly $m-$scalable,
\item[(ii)] There exists a subset $I\subset [M]$ with $\# I=m$ such that $0\in\co(F(\Phi_{I}))$,  respectively, $0\in ri\co(F(\Phi_{I}))$.
\item[(iii)] There exists a subset $I\subset [M]$ with $\# I = m$ for which there is no $h\in\R^d$ with $\ip{F(\varphi_k)}{h}> 0$ for all $k\in I$, respectively, with $\ip{F(\varphi_k)}{h}\geq 0$ for all $k\in I$, with at least one of the inequalities being strict.
\end{enumerate} 
\end{theorem}
 The details of the proof of this result can be found in  \cite[Theorem 3.9]{kop13c}. We point out however, that the equivalence of (i) and (ii) follows from considering $\co(F(\Phi_I))$ which is the polytope in $\R^N$ generated by the vectors $\{F(\varphi_k)\}_{k\in I}$. 
 
 By removing the normalization condition that the $\ell^1$ norm of the weights making a frame scalable is unity, Theorem~\ref{poly} can be stated in terms of the polyhedral cone $C(F(\Phi))$  generated by $F(\Phi)$.  This is the content of the following  result which was proved in \cite[Corollary 3.14]{kop13c}: 
\begin{cor} \cite[Corollary 3.14]{kop13c}\label{mainNd1}
Let $\Phi = \{\varphi_k\}_{k=1}^M\in \fr^*$, and let $N\leq m\leq M$ be fixed. Then the following conditions are equivalent:
\begin{enumerate}
\item[(i)]   $\Phi$ is strictly $m$-scalable .
\item[(ii)]  There exists $I \subset [M]$ with $\#I = m$ such that $C(F(\Phi_{I}))$ is not pointed.
\item[(iii)] There exists $I \subset [M]$ with $\#I = m$ such that $C(F(\Phi_{I}))^{\circ}$ is not blunt.
\item[(iv)]  There exists $I \subset [ M]$ with $\#I = m$ such that the interior of $C(F(\Phi_{I}))^{\circ}$ is  empty.
\end{enumerate}
\end{cor}

The map $F$ given in~\eqref{fctF}  is related to the diagram vector of \cite{hklw}, and was used in  \cite{cklmns12}  to give a different and equivalent characterization of  scalable frames which we now present. We start by presenting an interesting necessary condition for scalability in both $\R^N$ and $\C^N$ proved in \cite[Theorem 3.1]{cklmns12}:

\begin{theorem}\cite[Theorem 3.1]{cklmns12}\label{t:martalsuf}
Let $\Phi=\{\varphi_k\}_{k=1}^M \in \fr_u(M,N)$. If $\Phi \in \sca(M,N)$, then there is no unit vector $u\in \R^N$ such that $|\ip{u}{\varphi_k}|\geq \tfrac{1}{\sqrt{N}}$ for all $k=1, 2, \hdots, M$ and $|\ip{u}{\varphi_k}|> \tfrac{1}{\sqrt{N}}$ for at least one $k$. 
\end{theorem}

As pointed out in \cite{cklmns12} the condition in Theorem~\ref{t:martalsuf} is also sufficient only when $N=2$. We wish to compare this result to the following theorem that give a necessary and a (different) sufficient condition for scalability in $\R^N$, and these two conditions are necessary and sufficient only for $N=2$.

\begin{theorem}\cite[Theorem 4.1]{ckopw13}\label{pro:ns}
Let $\Phi\in\fr_{u}(M,N)$. Then the following hold:
\begin{enumerate}
\item[(a)] {\rm (}A necessary condition for scalability\,{\rm )} If $\Phi$ is scalable, then
\begin{equation}\label{e:nec}
\min_{\|d\|_2=1}\max_k |\ip{d}{\varphi_k}|\geq \frac{1}{\sqrt{N}}.
\end{equation} 
\item[(b)] {\rm (}A sufficient condition for scalability\,{\rm )} If
\begin{equation}\label{e:suff}
\min_{\|d\|_2=1}\max_k |\ip{d}{\varphi_k}|\geq \sqrt{\frac{N-1}{N}},
\end{equation}
then $\Phi$ is scalable.
\end{enumerate}
\end{theorem}

Clearly when $N=2$ the right hand sides of both~\eqref{e:nec} and~\eqref{e:suff} coincide leading to a necessary and sufficient condition.

Observe that~\eqref{e:nec} is equivalent to the fact that for each unit vector $d\in \R^N$, there exists $k =1, 2, \hdots, M$ such that $$|\ip{d}{\varphi_k}|\geq \frac{1}{\sqrt{N}}$$ which is different from the condition in theorem~\ref{t:martalsuf}.

We can now present the characterization of scalable  obtained in~\cite[Theorem 3.2]{cklmns12} and which is based on the Gramian of the diagram vectors. More precisely, for each $v \in \R^N$, we define the diagram vector\index{diagram vector} to be the vector $\tilde{v} \in \R^{N(N-1)}$ given by

\begin{equation}\label{diag-vec}
\tilde{v}= \tfrac{1}{\sqrt{N-1}}\begin{bmatrix}v(1)^2-v(2)^2\\ \vdots \\ v(N-1)^2-v(N)^2\\ \sqrt{2N}v(1)v(2)\\ \vdots \\ \sqrt{2N} v(N-1)v(N)\end{bmatrix},
\end{equation}
where the difference of the squares $v^2(i)-v^2(j)$ and the product $v(i)v(j)$ occur exactly once for $i<j, i=1, 2, \hdots, N-1.$ Using this notion, the following result was proved:

\begin{theorem}\cite[Theorem 3.2]{cklmns12} \label{mainthe-diag}
Let  $\Phi=\{\varphi_k\}_{k=1}M \in \fr_{u}$ be a frame of unit-norm vectors, and $\widetilde{G}=(\ip{\tilde{\varphi}_{k}}{\tilde{\varphi}_{\ell}})$ be the Gramian of the diagram vectors $\{\tilde{\varphi}_{k}\}_{k=1}^M$. Suppose that $\widetilde{G}$ is not invertible. Let $\{v_1, v_2, \hdots, v_\ell\}$ be a basis of the nullspace of $\widetilde{G}$ and set $$r_i:=\begin{bmatrix}v_1(i)\\ \vdots \\ v_\ell (i)\end{bmatrix},$$ for $i=1, 2, \hdots, M$. Then $\Phi$ is scalable if and only if $0 \not\in \co\{r_1, r_2, \hdots, r_M\}$. 

\end{theorem}

When a frame $\Phi=\{\varphi_k\}_{k=1}^M \subset \R^N$ is scalable, then  there exist $\{c_k\}_{k=1}^M \subset [0, \infty)$ such that $\{c_k \varphi_k\}_{k=1}^M$ is a tight frame. The nonnegative vector  $\omega=\{c_k^2\}_{k=1}^M$ is called a \emph{scaling} of $\Phi$ \cite{cc12}. The scaling $\omega=\{c_k^2\}_{k=1}^M \subset [0, \infty)$ is said to be a \emph{minimal scaling}\index{minimal scaling} if $\{\varphi_k: c_k^2>0\}$ has no proper subset which is scalable. The notion of minimal scalings has recently found some very interesting applications on  some structural decomposition of frames; see,  \cite[Section 4]{ccnst14} for more details.  It turns out that finding the scalings of a scalable frame can be reduced to finding its minimal scalings. More specifically, the following result was proved in \cite[theorem 3.5]{cc12}:

\begin{theorem}\cite[Theorem 3.5]{cc12}
Suppose $\Phi=\{\varphi_k\}_{k=1}^M \in \mathcal{F}_u$ is a scalable frame, and let $\omega=\{\omega_k\}_{k=1}^M \subset [0, \infty)$ be one of its  minimal scalings. Then $\{\varphi_k \varphi_k^T: \omega_k>0\}$ is linearly independent. Furthermore, every scaling of $\Phi$ is a convex combination of minimal scalings. 
\end{theorem}

\subsection{A geometric condition for scalability}\label{subsec2.3} The two dimensional case we examined earlier (Example~\ref{example2-3} and Example~\ref{example2-4}) indicates that a frame is not scalable 
when the frame vectors ``cluster'' in certain ``small''   plane regions. In fact, broadly speaking, the frame is not scalable if its  vectors lies in a double cone $C\cap (-C)$ with a ``small'' aperture. This was  formalized in theorem~\ref{poly} and Corollary~\ref{mainNd1}. We can further exploit these results to give a more formal geometric characterization of scalable frames. 

To begin, we rewrite  (iii) of Theorem~\ref{poly} in the following form. For $x=(x_k)_{k=1}^{N}\in \R^N$ and $h=(h_k)_{k=1}^d \in \R^d$, we have that
\begin{equation}\label{poly2}
\ip{F(x)}{h} = \sum_{\ell=2}^{N}h_{\ell-1}(x_1^2 - x_{\ell}^2) + \sum_{k=1}^{N-1}\sum_{\ell=k+1}^{N}h_{k(N - 1 - (k-1)/2) +\ell - 
1}x_{k}x_{\ell}.
\end{equation}
Consequently, fixing $h \in \R^d$, $\ip{F(x)}{h}$  is  a homogeneous polynomial of degree $2$ in $x_1, x_2, \hdots, x_N$. Denote the set of all polynomials of this form by $\boldsymbol{P}_{2}^N$. Then   $\boldsymbol{P}_{2}^N$ can be identified with  the  subspace of real symmetric $N\times N$ matrices whose trace is $0$. Indeed, for each $N\geq 2$, and each $p\in\boldsymbol{P}_{2}^N$,
$$
p(x) = \sum_{\ell=2}^{N}a_{\ell-1}(x_1^2 - x_{\ell}^2) + \sum_{k=1}^{N-1}\sum_{\ell=k+1}^{N}a_{k(N - (k+1)/2) +\ell - 1}x_{k}x_{\ell},
$$
we have $p(x) = \ip{Q_px}{x}$, where $Q_p$ is the symmetric $N\times N$-matrix with entries
$$
Q_{p}(1,1) = \sum_{k=1}^{N-1}a_k, \qquad Q_{p}(\ell,\ell) = -a_{\ell-1}\quad\text{for }\ell = 2,3,\ldots,N
$$
and
$$
Q_{p}(k,\ell) = \frac 1 2 a_{k(N - (k+1)/2) +\ell - 1}\quad\text{for }k = 1,\ldots,N-1,\;\ell = k+1,\ldots,N.
$$

Thus, $\ip{F(x)}{h}=\ip{Q_{h}x}{x}=0$ defines a quadratic surface\index{quadratic surface} in $\R^N$, and condition (iii) in Theorem~\ref{poly} stipulates that for $\Phi$ to be scalable, one cannot find such a quadratic surface such that the frame vectors (with index in $I$) all lie  on (only) ``one side'' of this surface. By taking the contrapositive statement  we arrived at the following result that was proved differently in  \cite[Theorem 3.6]{kopt12}. In particular, it provides a characterization of non-scalability of  finite frames, and  we shall use it to give  a very interesting geometric condition on the frame vectors for non-scalable frames. 

\begin{theorem}\cite[Theorem 3.6]{kopt12}\label{t:quadric}
Let $\Phi = \{\varphi_k\}_{k=1}^M\in \fr^*.$ Then the following statements are equivalent.
\begin{enumerate}
\item[{\rm (i)}]   $\Phi$ is not scalable.
\item[{\rm (ii)}]  There exists a symmetric matrix $Y\in\R^{N\times N}$ with $\tr(Y) < 0$ such that $\ip{Y \varphi_k}{\varphi_k}\ge 0$ for all $k= 1,\ldots,M$.
\item[{\rm (iii)}] There exists a symmetric matrix $Y\in\R^{N\times N}$ with $\tr(Y) = 0$ such that $\ip{Y\varphi_k}{\varphi_k } >0$ for all $k = 1,\ldots,M$.
\end{enumerate}
\end{theorem}

To derive the geometric condition for non-scalability we need some set up. It is not difficult to see that each symmetric $N\times N$ matrix $Y$  in (iii) of theorem~\ref{t:quadric} corresponds to a quadratic surface. We call this surface \emph{a conical zero-trace quadric}. The exact definition of such quadratic surface is 

\begin{definition}\cite[Definition 3.4]{kopt12}
Let the \emph{class of conical zero-trace quadrics} $\calC_N$ be defined as the family of sets
\begin{equation}\label{eq:quadricvariety}
\left\{x\in\R^N : \sum_{k=1}^{N-1} a_k \ip{x}{e_k}^2 = \ip{x}{e_N}^2\right\},
\end{equation}
where $\{e_k\}_{k=1}^N$ runs through all orthonormal bases of $\R^N$ and $(a_k)_{k=1}^{N-1}$
runs through all tuples of elements in $\R\setminus\{0\}$ with $\sum_{k=1}^{N-1} a_k = 1$.
\end{definition}

We define the {\em interior} of the conical zero-trace quadric in \eqref{eq:quadricvariety}, by 
\[
\left\{x\in\R^N : \sum_{k=1}^{N-1} a_k \ip{x}{e_k}^2 < \ip{x}{e_N}^2\right\},
\]

and the {\em exterior} of the conical zero-trace quadric in \eqref{eq:quadricvariety} by 

\[  \left\{x\in\R^N : \sum_{k=1}^{N-1} a_k \ip{x}{e_k}^2 > \ip{x}{e_N}^2\right\}. 
\]

It is then easy to see that Theorem~\ref{t:quadric} is equivalent to the following result established  in \cite[theorem 3.6] {kopt12}

\begin{theorem}\cite[Theorem 3.6] {kopt12}\label{t:geometry}
Let $\Phi\subset\R^N\setminus\{0\}$ be a frame for $\R^N$. Then the following conditions are equivalent.
\begin{itemize}
\item[{\rm (i)}]   $\Phi$ is not scalable.
\item[{\rm (ii)}]  All frame vectors of $\Phi$ are contained in the interior of a conical zero-trace quadric of $\calC_N$.
\item[{\rm (iii)}]  All frame vectors of $\Phi$ are contained in the exterior of a conical zero-trace quadric of $\calC_N$.
\end{itemize}
\end{theorem}

\begin{figure}[h]
\includegraphics[width=5cm]{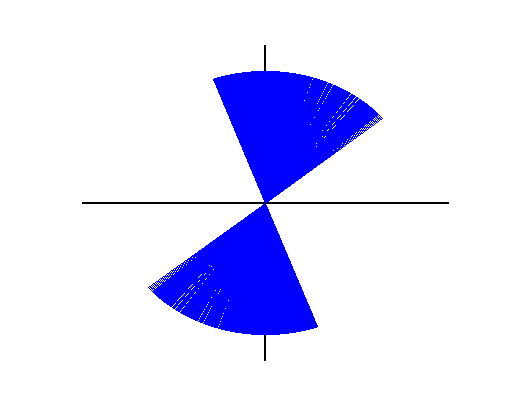}
\hspace*{0.5cm}
\includegraphics[width=5cm]{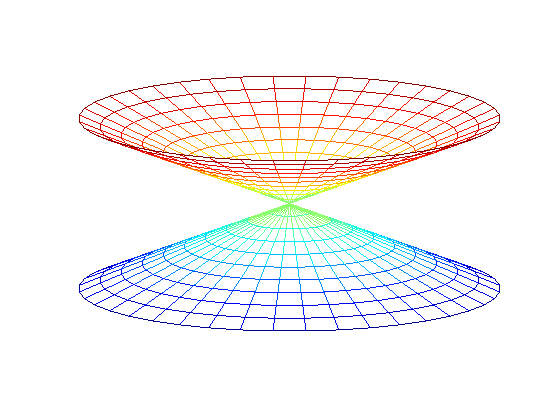}
\put(-250,-15){(a)}
\put(-80,-15){(b)}
\caption{(a) shows a sample region of vectors of a non-scalable frame in $\R^2$. (b) 
shows example of $\mathcal{C}^-_3$ and $\mathcal{C}^+_3$ which determine sample regions in $\R^3$.}
\label{fig-5}
\end{figure}

The geometric meaning of this result is best illustrated by considering frames in $\R^2$ and $\R^3$, in which case  the sets $\mathcal{C}_N$ for $N=2, 3$, have  very simple descriptions given in \cite{kopt12}. For our purposes here it suffices to say that  each set in $\mathcal{C}_2$  is the boundary surface of a {\it quadrant cone} in $\R^2$, i.e., the union of two orthogonal one-dimensional subspaces (lines through the origin) in $\R^2$. The sets in $\mathcal{C}_3$ are the boundary surfaces of a particular class of {\it elliptical cones} in $\R^3$. We give examples of sets in $\mathcal{C}_N$ $N=2,3$ in Figure~\ref{fig-5} (a) and (b).

We can now state the following corollary that give a clear geometric insight into the set of non-scalable frames. In particular, the frame vectors cannot lie in a ``small cone''. 

\begin{cor}\cite[Corollary 3.8]{kopt12}\label{cone2-3}
\begin{itemize}
\item[{\rm (i)}] A frame $\Phi\subset\R^2\setminus\{0\}$ for $\R^2$ is not scalable if and only if there exists an open quadrant cone
which contains all frame vectors of $\Phi$.
\item[{\rm (ii)}] A frame $\Phi\subset\R^3\setminus\{0\}$ for $\R^3$ is not scalable if and only if all frame vectors of $\Phi$ are
contained in the interior of an elliptical conical surface with vertex $0$ and intersecting the corners of a rotated unit cube.
\end{itemize}
\end{cor}

\subsection{Scalable frames and Fritz John theorem}\label{subsec2.4}
The last characterization of scalable frames we should discuss is based on Fritz John's ellipsoid theorem\index{Fritz John}. Before we state this theorem, we recall from Section~\ref{subsec2.2} that given a set of points  $Y =\{y_k\}_{k=1}^{L}\subset \R^{N}$,  $P_{Y}$ is the polytope generated by $Y$. 

Given an $N\times N$ positive definite matrix $X$ and a point $c \in \R^N$, we define an $N$-dimensional ellipsoid\index{ellipsoid} centered at $c $ as 
$$E(X,c)=c+X^{-1/2}(B)=\{v: \ip{ X(v-c)}{ (v-c)}\leq 1\},$$ 
where $B$ is the unit ball in $\R^N$.
We recall that the volume of the ellipsoid is given by  $\text{Volume}(E(X,c))=(\det(X))^{-1/2}\omega_N$, where $\omega_N$ is the volume of the unit ball in $\R^N$.

A convex body $K \subset \R^N$ is a nonempty compact convex subset of $\R^N$. It is well-known that for any convex body $K$ with nonempty interior in $\R^N$ there is a unique ellipsoid  of minimal volume\index{ellipsoid!ellipsoid of minimal volume} containing $K$; e.g., see \cite[Chapter 3]{ntj89}. We refer to \cite{ kb92, kb97,  og10, fj48,  ntj89} for more on these extremal ellipsoids\index{ellipsoid!Fritz John ellipsoid}. Fritz John ellipsoid theorem \cite{fj48} gives a description of this ellipsoid. More specifically:

\begin{theorem}\cite[Section 4]{fj48}\label{FJohn}  Let $K \subset B=B_2^N(0,1)$ (unit ball in $\R^N$) be a convex  body with nonempty interior. Then $B$ is the ellipsoid of minimal volume containing $K$ if and only if there exist $\{\lambda_k\}_{k=1}^{m} \subset (0, \infty)$ and $\{u_{k}\}_{k=1}^{m} \subset \partial{K}\cap S^{N-1}$, $m\geq N+1$ such that 
\begin{enumerate}
\item[(i)] $\sum_{k=1}^{m}\lambda_k u_k =0$
\item[(ii)] $x=\sum_{k=1}^{m}\lambda_{k} \ip{x}{u_k}u_k, \forall x \in \R^N$
\end{enumerate} where $\partial K$ is the boundary of $K$ and $S^{N-1}$ is the unit sphere in $\R^N$. In particular, the points $u_k$ are contact points of $K$ and $S^{N-1}$. 
\end{theorem}

Observe that (ii) of Theorem~\ref{FJohn} can be written as $$x=\sum_{k=1}^m \ip{x}{\sqrt{\lambda_k}u_k}\sqrt{\lambda_k}u_k\quad {\textrm for \, all} \quad  x \in \R^N$$ which is equivalently to saying that the vectors $\{u_k\}_{k=1}^m$ form a scalable frame. The difficulty in applying this theorem lies in the fact that determining the contact points $u_k$ and the ``multipliers'' $\lambda_k$ is an extremely difficult problem. Nonetheless, we can apply this result to our question since we  consider the convex body generated by the frame vectors, in which case these contact points are subset of the frame vectors. 
In particular, to apply the Fritz John theorem to the scalability problem, we consider   a frame $\Phi\in\fr_u(M,N)$ of $\R^N$ consisting of unit norm vectors. We define the associated symmetrized frame as \[
\Phi_{\rm Sym} := \{\varphi_k\}_{k=1}^M\cup\{-\varphi_k\}_{k=1}^M,
\] and we denote the ellipsoid of minimal volume  circumscribing the convex hull of the symmetrized frame  $\Phi_{\rm Sym}$  by $E_{\Phi}$ and refer to it as the \emph{minimal ellipsoid\index{ellipsoid!minimal ellipsoid}  of $\Phi$}. Its `normalized' volume is defined by
\[
V_{\Phi} := \frac{\Vol(E_{\Phi})}{\omega_N}. 
\]
Clearly,  $V_\Phi \le 1$, and it is shown in \cite[theorem 2.11]{ckopw13} that equality holds if and only if the frame is scalable. That is, we have

\begin{theorem}\cite[Theorem 2.11]{ckopw13}\label{pro:iff}
A frame $\Phi\in\fr_u(M,N)$ is scalable if and only if its minimal ellipsoid is the $N$-dimensional unit ball, in which case $V_\Phi = 1$.
\end{theorem}

\begin{remark}
Given a unit-norm frame $\Phi$, the number $V_\Phi$ defined above is one of a few measures of scalability introduced in \cite{ckopw13}. These are numbers that measure how close to being scalable a frame $\Phi$ is. For example, if for a given $\Phi$, $V_\Phi <1$, then the farther away from $1$ it is, the less scalable is $\Phi$. Thus $V_\Phi$ along with these other measure of scalability can be used to define ``almost'' scalable frames. We refer to \cite{ckopw13} for details. 
\end{remark}

Using the geometric characterization of scalable frames by $V_\Phi$ one can  define the following  equivalence relation on $\fr_u(M,N)$:  $\Phi,\Psi\in\fr_u(M,N)$ are  equivalent if and only if $V_{\Phi} = V_{\Psi}$. We denote each equivalence class by the unique volume for all its members. Specifically, for any $0<a\leq 1$, the class $P[M, N, a]$ consists of all $\Phi\in \fr_u(M,N)$ with $V_{\Phi}=a$. Then, $\sca(M,N)=P[M,N,1]$. This also allows a parametrization of $\fr_u(M,N)$:
$$
\fr_u(M,N) = \bigcup_{a\in (0,1]}P[M,N,a].$$

\section{Probabilistic frames}\label{sec3}

Definition~\ref{framedef} introduces   frames  from a linear algebra perspective  through  their spanning properties. However,  frames can also be viewed as  point masses  distributed  in $\R^d$. In this section we survey a measure theoretical, or more precisely a probabilistic,  description  of frames. In particular, in Section~\ref{subsec3.1} we define probabilistic frames and collect some of their elementary properties. In Section~\ref{subsec3.2} We define the probabilistic frame potential and investigate its minimizers. We generalize this notion in Section~\ref{subsec3.3} to the concept of $p^{th}$ frame potentials and discuss their minimizers. Probabilistic analogs of these potentials are considered in Section~\ref{subsec3.4}.   This section can be considered as a companion to \cite{ehloko11} where many of the results we stated below first appeared. 

\subsection{Definition and elementary properties}\label{subsec3.1} Before defining probabilistic frames we first collect some definitions needed in the sequel.

Let $\mathcal{P}:=\mathcal{P}(\mathcal{B},\R^N)$ denote the collection of probability measures on $\R^N$ with respect to the Borel $\sigma$-algebra $\mathcal{B}$.  
Let $$
\mathcal{P}_{2}:= \mathcal{P}_{2}(\R^{N})=\bigg\{ \mu \in \mathcal{P}: M_{2}^{2}(\mu):=\int_{\R^{N}}\nm{x}{}^{2}d\mu(x) < \infty\bigg\}$$ 
be the set of all probability measures with finite second moments.  Given $\mu, \nu \in \mathcal{P}_2$, let 
$\Gamma(\mu, \nu)$ be the set of all Borel probability measures  $\gamma$ on $\R^N \times \R^N$ whose marginals are $\mu$ and $\nu$, respectively, i.e., $\gamma(A\times \R^N)=\mu(A)$ and $\gamma(\R^N \times B) = \nu(B)$ for all Borel subset $A, B$ in $\R^N$. The space 
$\mathcal{P}_{2}$ is equipped with the $2$-\emph{Wasserstein metric}  given by 
\begin{equation}\label{wmetric}
W_{2}^{2}(\mu, \nu):=\min\bigg\{\int_{\R^N \times \R^N}\nm{x-y}{}^{2}d\gamma(x, y), \gamma \in \Gamma(\mu, \nu)\bigg\}.
\end{equation}
 It is known that the minimum defined by~\eqref{wmetric} is attained at a  measure $\gamma_0 \in \Gamma(\mu, \nu)$, that is:
$$W_{2}^{2}(\mu, \nu)=\int_{\R^N \times \R^N}\nm{x-y}{}^{2}d\gamma_0(x, y). $$ 
  We refer to \cite[Chapter 7]{ags05}, and  \cite[Chapter 6]{vill09} for more details on the Wasserstein spaces.  

\begin{definition}\label{defpf}
A Borel probability measure $\mu \in \mathcal{P}$  is a \emph{probabilistic frame} if there exist $0<A\leq B < \infty$ such that for all $x\in\R^N$ we have 
\begin{equation}\label{pfineq}
 A\|x\|^2 \leq \int_{\R^{N}} |\langle x,y\rangle |^2 d\mu (y) \leq B\|x\|^2. 
 \end{equation}
The constants $A$ and $B$ are called \emph{lower and upper probabilistic frame bounds}, respectively.
When $A=B,$ $\mu$ is called a \emph{tight probabilistic frame}.   
\end{definition}

 It follows from Definition~\ref{defpf}  that the upper inequality in~\eqref{pfineq} holds if and only if $\mu \in \mathcal{P}_2$. With a little more work one shows that  the lower inequality holds whenever the linear span of the support of the probability measure $\mu$ is $\R^N$.

Assume that  $\mu$ is a tight probabilistic frame, in which case equality holds in~\eqref{pfineq}. Hence, choosing $x=e_k$ where $\{e_k\}_{k=1}^N$ is the standard orthonormal basis for $\R^N$  leads to $$A\|e_k\|^2=A=\int_{\R^N}\ip{e_k}{y}^2 d\mu(y).$$  Therefore, $$NA=\sum_{k=1}^N A\|e_k\|^2= \int_{\R^N}\sum_{k=1}^N \ip{e_k}{y}^2 d\mu(y)=\int_{\R^N}\|y\|^2d\mu(y)=M_2(\mu)^2.$$ Consequently,   for a  tight probabilistic frame $\mu$, $A=\tfrac{M_{2}(\mu)^2}{N}$. 

  These observations can are  summarized in the following result whose proof can be found in \cite{ehloko11}

\begin{theorem}\cite[Theorem 12.1]{ehloko11}\label{thm:charapf}
A Borel probability measure $\mu \in \mathcal{P}$  is a \emph{probabilistic frame}\index{probabilistic frame} if and only if $\mu \in \mathcal{P}_{2}$ and $E_{\mu}=\R^N$, where  $E_{\mu}$ denotes the linear span of $\supp(\mu)$ in $\R^N$. Moreover, if $\mu$ is a tight probabilistic frame\index{probabilistic frame!tight probabilistic frame}, then the frame bound is given by $$A=\tfrac{1}{N}M_{2}^{2}(\mu)=\tfrac{1}{N}\int_{\R^{N}}\|y\|^{2}d\mu(y).$$
\end{theorem}

We now consider some examples of probabilistic frames.

\begin{example}\label{examp-probframes}
\begin{enumerate}
\item[(a)] A set  $\Phi=\{\varphi_k\}_{k=1}^{M} \subset \R^N$ is a frame   if and only if the probability measure $\mu_{\Phi}=\tfrac{1}{M}\sum_{k=1}^{M}\delta_{\varphi_k}$ supported by the set $\Phi$ is a probabilistic frame, where $\delta_\varphi$ denotes the Dirac measure supported at $\varphi\in\R^N$. 
\item[(b)] More generally, let $a=\{a_k\}_{k=1}^M \subset (0, \infty)$ with $\sum_{k=1}^Ma_k=1$. A set  $\Phi=\{\varphi_k\}_{k=1}^{M} \subset \R^N$ is a frame   if and only if the probability measure $\mu_{\Phi, a}=\sum_{k=1}^{M}a_k\delta_{\varphi_k}$ supported by the set $\Phi$ is a probabilistic frame.
\item[(c)] By symmetry consideration one also shows that the uniform distribution on the unit sphere  $S^{N-1}$ in $\R^N$  is a tight probabilistic frame \cite[Proposition 3.13]{ehler10}. That is, denoting the probability measure on $S^{N-1}$ by $d\sigma$ we have that for all $x \in \R^N$, $$\tfrac{\|x\|^2}{N}=\int_{\R^N}\ip{x}{y}^2d\sigma(y).$$
\end{enumerate}
\end{example}

In the framework of the  Wasserstein metric, many properties of probabilistic can be proved. For example, if we denote by $P(A, B)$ the set of probabilistic frames with frame bounds $0<A \leq B < \infty$, then the following result was proved in \cite[Proposition 1]{cwko15}:

\begin{proposition}\label{pfab} \cite[Proposition 1]{cwko15}
$P(A,B)$ is a nonempty, convex, closed subset of $\mathcal{P}_2(\R^N)$.
\end{proposition}
Other results including a probabilistic treatment of the frame scalability problem also appeared in \cite{cwko15}. Furthermore, in \cite{cw14} an optimal transport approach to minimizing a frame potential\index{frame potential} that generalizes the Benedetto and Fickus potential was developed. In the process, the smoothness (in the Wasserstein metric) of this potential was derived.

Probabilistic frames   can be analyzed in terms of a corresponding analysis operator and its adjoint the synthesis operator.
Indeed, let $\mu \in \mathcal{P}$ be a probability measure. The \emph{probabilistic analysis operator} is given by
 $$
 T_\mu: \R^N \rightarrow L^2(\R^N,\mu), \quad x\mapsto \ip{ x}{\cdot}.$$
Its adjoint operator is defined by  
 $$
  T^*_\mu: L^2(\R^N,\mu) \rightarrow\R^N , \quad f\mapsto \int_{ \R^N } f(x)xd\mu(x)$$ and is called 
the \emph{probabilistic synthesis operator}, where the above  integral is vector-valued.  The \emph{probabilistic frame operator}\index{frame operator!probabilistic frame operator} of $\mu$ is $$S_\mu=T^*_\mu T_\mu,$$ and one easily verifies that 
\begin{equation*}
S_\mu:\R^N\rightarrow \R^N,\qquad S_\mu (x) = \int_{\R^N} \ip{ x}{y} y d\mu(y).
\end{equation*}

If $\{e_j\}_{j=1}^N$ is the canonical orthonormal basis for $\R^N$, then 
 $$
S_\mu e_i =  \sum_{j=1}^N m_{i,j}(\mu) e_j, $$ where $$m_{i,j}= \int_{\R^N}  y^{(i)} y^{(j)} d\mu(y)$$ is the $(i, j)$ entry of the matrix of second moments of $\mu$. Thus, the probabilistic frame operator is the matrix of second moments of $\mu$. Consequently, the following results proved in \cite{ehloko11} follows. 

\begin{proposition}\cite[Proposition 12.4]{ehloko11}\label{pf:operator}
Let  $\mu\in\mathcal{P}$, then 
$S_\mu$ is well-defined (and hence bounded) if and only if $$
 M_2(\mu)<\infty.$$ Furthermore, 
$\mu$ is a probabilistic frame if and only if $S_\mu$ is  positive definite.  
\end{proposition}

If $\mu$ is a probabilistic frame then $S_\mu$ is invertible. Let $\tilde{\mu}$ be the \emph{push-forward} of $\mu$ through $S_{\mu}^{-1}$ given by  $$\tilde{\mu}=\mu \circ S_\mu.$$  In particular, given any Borel set $B\subset \R^N$ we have  $$\tilde{\mu}(B) = \mu((S_{\mu}^{-1})^{-1} B) = \mu(S_{\mu}B).$$     Equivalently, $\tilde{\mu}$ can be defined via integration. Indeed,  if $f$ is a continuous bounded function on $\R^N$, $$\int_{\R^N} f(y)d\tilde{\mu}(y) = \int_{\R^N} f(S_{\mu}^{-1}y) d\mu(y).$$
In fact,  $\tilde\mu$ is also a probabilistic frame (with bounds $1/B\leq 1/A$) called  the \emph{probabilistic canonical dual frame}\index{probabilistic frame!probabilistic canonical dual frame} of $\mu$. Similarly, when $\mu$ is a probabilistic frame, $S_\mu$ is positive definite, and its square root exists. The push-forward of $\mu$ through $S_{\mu}^{-1/2}$ is given by $$\mu^{\dagger} (B) =\mu(S^{1/2} B)$$ for each Borel set in $\R^N$. The properties of these probability measures are summarized in the following result. We refer to \cite[Proposition 12.4]{ehloko11} and \cite[Proposition 12.5]{ehloko11} for details. 

\begin{proposition}\label{can-dual-reconstruction}
Let  $\mu\in\mathcal{P}$ be a probabilistic frame with bounds $0<A\leq B < \infty.$ Then:
\begin{enumerate}
\item[(a)] $\tilde\mu$ is a probabilistic frame with frame bounds $1/B\leq 1/A$.
\item[(b)] $\mu^{\dagger} $ is a tight probabilistic frame. 
\end{enumerate}
Consequently,  for each $x \in \R^N$ we have: 
\begin{equation}\label{reconsppf}
\int_{\R^N} \langle x, y\rangle \, S_{\mu}y \, d\tilde{\mu}(y)= \int_{\R^N} \langle  S_{\mu}^{-1} x, y \rangle \, y \, d\mu(y) =S_{\mu}S_{\mu}^{-1}(x)=x, 
\end{equation} 
and 
\begin{equation}\label{reconsptf}
\int_{\R^N} \langle x, y\rangle \, y \, d\mu^{\dagger} (y)= \int_{\R^N} \langle  S_{\mu}^{-1/2} x, y \rangle \, S_{\mu}^{-1/2}y \, d\mu(y) =S_{\mu}^{-1/2}S_{\mu} S_{\mu}^{-1/2}(x)=x.
\end{equation} 
\end{proposition}

It is worth noticing that ~\eqref{reconsppf} is the analog of the frame reconstruction formula~\eqref{framerecons} while 
~\eqref{reconsptf} is analog of~\eqref{tfrecon2}.

In the context of probabilistic frames, the \emph{probabilistic Gram operator}\index{Gram operator!probabilistic Gram operator}, or the \emph{probabilistic Gramian} of $\mu$, is  the compact integral operator defined on $L^2 (\R^N, \mu)$ by
\begin{equation*}G_{\mu}f(x)=T_{\mu}T^{*}_{\mu}f(x)=\int_{\R^{N}}K(x,y)f(y)d\mu(y)=\int_{\R^{N}}\langle x,y\rangle f(y)d\mu(y).
\end{equation*}  It is immediately seen that 
 $G_\mu$  is an integral operator with kernel given by $K(x, y)=\ip{x}{y}$, which is continuous and in $L^{2}(\R^{N}\times \R^{N}, \mu\otimes \mu ) \subset 
L^{1}(\R^{N}\times \R^{N}, \mu\otimes \mu )$, where $\mu\otimes \mu$  is the product measure of $\mu$ with itself. Consequently, $G_\mu$ is a  trace class and Hilbert-Schmidt operator. Moreover, for any $f \in L^{2}(\R^{N}, \mu)$, $G_{\mu}f$ is a uniformly continuous function on $\R^N$. As well-known,  $G_\mu$ and $S_\mu$ have a common spectrum except for the $0$. In fact, in the next proposition we collect the properties of $G_\mu$:

\begin{proposition}\cite[Proposition 12.4]{ehloko11}\label{gram-proposition}
Let $\mu \in \mathcal{P}$ then $G_\mu$ is a  trace class and Hilbert-Schmidt operator on $L^2(\R^N)$. The eigenspace corresponding to the eigenvalue $0$ has infinite dimension and consists of all functions  $0\neq f \in L^2(\R^N, \mu)$ such that  $$\int_{\R^N}y f(y) d\mu(y)=0.$$ 

\end{proposition}

While new finite frames can be generated from old ones via (linear) algebraic operations, the setting of probabilistic frames allows one to use analytical tools to construct new probability frames from old ones. For example, it was shown in \cite{ehloko11} when the convolution of a probabilistic frame and a probability measure yields a probabilistic frames. The following is a summary of some of the results proved in \cite{ehloko11}. 

\begin{proposition}\cite[Theorem 2 $\&$ Proposition 2]{ehloko11} The following statements hold:
\begin{enumerate}
\item[(a)] Let $\mu\in \mathcal{P}_{2}$ be a probabilistic frame and let $\nu \in \mathcal{P}_{2}$. If $\supp(\mu)$ contains at least $N+1$ distinct vectors, then $\mu\ast\nu$ \index{probabilistic frame!convolution of probabilistic frames} is a probabilistic frame.
\item[(b)] Let $\mu$ and $\nu$ be tight probabilistic frames. If $\int_{\R^{N}}y d\nu(y)=0$, then $\mu \ast \nu$ \index{probabilistic frame!convolution of probabilistic frames} is also a tight probabilistic frame. 
\end{enumerate}
\end{proposition}

\subsection{Probabilistic frame potential}\label{subsec3.2} One of the motivations of probabilistic frames lies in Benedetto and Fickus's characterization of the FUNTFs as  the minimizers of the  frame potential~\eqref{potineq}. In describing their results, they motivated it from a physical point of view drawing a parallel to Coulomb's law. It was then clear that the notion of frame potential carries significant information about frames, and can be viewed as describing the interaction of the frame vectors under some ``physical force.''   This in turn partially motivated the introduction of a probabilistic analog to the frame potential in \cite{ehloko10}. Furthermore, the probabilistic frame potential\index{frame potential!probabilistic frame potential} that we introduce below, can be viewed in the framework of other potential functions, e.g., those investigated by Bj\"orck in \cite{Bjoerck:1955aa}.  In this section we review the properties of the probabilistic frame potential investigating in particular its minimizers. The framework of the Wasserstein metric space $(\mathcal{P}_2, W_2)$ also offers the ideal setting to investigate this potential and certain of its generalizations. While we should not report of this analysis here, we shall nevertheless introduce certain  generalizations of the probabilistic frame potential whose minimizers are better understood in the context of the Wasserstein metric spaces. 

But we first start with the definition of the probabilistic frame potential. 

\begin{definition}\label{probfpot}
The \emph{probabilistic frame potential} is the nonnegative function defined on $\mathcal{P}$ and given by 
\begin{equation}\label{eqpfp}
\pfp(\mu)=\iint_{\R^{N}\times \R^{N}}|\ip{ x}{y} |^{2}\, d\mu(x)\, d\mu(y),
\end{equation} for each $\mu \in \mathcal{P}$. 
\end{definition}

The following proposition is an immediate consequence of the above definition: 

\begin{proposition}\label{tracegramian}
Let $\mu \in \mathcal{P}$, then $\pfp(\mu)$ is the Hilbert-Schmidt norm of the probabilistic Gramian operator $G_\mu$, that is 
$$\nm{G_\mu}_{HS}^{2}=\iint_{\R^{d}\times \R^{d}} \ip{x}{y}^{2}d\mu(x) d\mu(y).$$

Furthermore, if $\mu \in \mathcal{P}_2,$ (which is the case when $\mu$ is a probabilistic frame) then we have $$\pfp(\mu) \leq M_{2}^{4}(\mu)< \infty.$$ 
\end{proposition}

We recall from Definition~\ref{defpf} that $\mu$ is a tight probabilistic frame if $$\int_{\R^N}\ip{x}{y}^2d\mu(y)=\tfrac{M_2(\mu)^2}{N}\|x\|^2$$ for all $x\in \R^N$. Integrating this equation with respect to $x$ leads to $$\iint_{\R^N \times \R^N}\ip{x}{y}^2d\mu(x)\, d\mu(y)=\pfp(\mu)=\tfrac{M_2^4(\mu)}{N}.$$ It turns out that this value is the absolute lower bound to the probabilistic frame potential. 

\begin{theorem}\cite[Theorem 3]{ehloko11}\label{minfp} 
Let $\mu \in \mathcal{P}_{2}$ be such that $M_{2}(\mu)=1$  and set  $E_{\mu}=span(\supp(\mu))$, then
the following estimate holds
\begin{equation}\label{estpfp1}
  \pfp(\mu) \geq 1/n
 \end{equation}  
 where $n$ is the number of nonzero eigenvalues of $S_\mu$. Moreover,  equality holds if and only if $\mu$ is a tight probabilistic frame for $E_{\mu}$.
 
 In particular, given any probabilistic frame $\mu \in \mathcal{P}_{2}$ with  $M_{2}(\mu)=1$, we have $$\pfp(\mu)\geq 1/N$$ and equality holds if and only if $\mu$ is a tight probabilistic frame.
  \end{theorem}
  
  The proof of this result can be found in \cite[Theorem 3]{ehloko11}. Recently, a very simple and elementary proof of the last part of the result was given in \cite[Theorem 5]{ckko15}. 
  Furthermore, in \cite{cw14} an optimal transport approach to minimizing a modification of the probabilistic frame potential was considered and showed great promise to analyze other potential functions in frame theory. Moreover, this approach has a natural numerical part that could be used as a gradient descent-type method to numerically find the minimizers of the $\pfp$ and its generalization.

\subsection{The $p^{th}$  frame potentials}\label{subsec3.3}  The techniques used to prove Theorem~\ref{minfp}  can be used to investigate the minimizers of other related potential functions, especially when there are defined for probability measures supported on compact sets, such as on the unit sphere. In this section, we define a family of (deterministic) potentials and describe  their minimizers. The probabilistic analogs  of these results will follow in the next section.  

To motivate our definition, we recall the following result due to Strohmer and Heath \cite{strhea}, and we refer to \cite{welch74} for historical perspectives on this result. 

\begin{theorem}\cite[Theorem 2.3]{strhea}
For any frame $\Phi=\{\varphi_k\}_{k=1}^M \subset \R^N$  with $\|\varphi_k\|=1$, we have 
\begin{equation}\label{welcbound}
\max_{k\neq \ell}|\ip{\varphi_k}{\varphi_\ell}|\geq \sqrt{\tfrac{M-N}{N(M-1)}},
\end{equation}
 and equality hold if and only if $\Phi$ is a FUNTF such that 
 \begin{equation}\label{equiangular}
 |\ip{\varphi_k}{\varphi_\ell}|=\sqrt{\tfrac{M-N}{N(M-1)}}
 \end{equation} when $k\neq \ell$.  Furthermore, equality can hold only when $M\leq \tfrac{N(N+1)}{2}.$
\end{theorem}
A FUNTF that satisfies~\eqref{equiangular}  is termed an \emph{ equiangular tight frame (ETF)}\index{frame!equiangular tight frame}. Note that the left-hand side of~\eqref{welcbound} can be viewed as a potential function of $\Phi$. Indeed, this is the so-called \emph{coherence}\index{coherence} of $\Phi$ that we shall defined for reasons to be evident later as 

\begin{equation}\label{coherence}
\FP_{\infty, M}(\Phi)=\max_{ k\neq \ell} |\ip{\varphi_{k}}{\varphi_{\ell}}|
\end{equation} for $\Phi=\{\varphi_k\}_{k=1}^M \subset S^{N-1}$. In fact, $\FP_{\infty, N}(\Phi)$ as well as the frame potential $FP$ given in~\eqref{potineq} are members of the family of the \emph{ $p^{th}$ frame potentials} defined by:

\begin{definition}\label{pframepotential}
Let $M$ be a positive integer, and $0<p< \infty$. Given a collection of unit vectors $\Phi=\{\varphi_k\}_{k=1}^{M}\subset S^{N-1}$, the \emph{$p$-frame potential}\index{$p$-frame potential} is the functional
  \begin{equation}\label{eq:pth frame potential}
  \FP_{p, M}(\Phi)=\sum_{k, \ell=1}^{M}|\ip{\varphi_k}{\varphi_\ell} |^p. 
 \end{equation}
When, $p=\infty$, the definition reduces to

\begin{equation*}
\FP_{\infty, M}(\Phi )=\max_{k\neq \ell} |\ip{\varphi_k}{\varphi_\ell} |.
\end{equation*}
\end{definition}

It is clear that $\FP_{p, M}$ and its minimizers are also functions of $N$, the dimension of the underlying space. However, to keep the notations simple, we shall not make explicit the dependence on $N$, unless it is necessary. 

The case $p=2$ corresponds to the frame potential  $FP$ given in~\eqref{potineq}. As mentioned above, $\FP_{\infty, M}(\Phi)$ is the coherence of $\Phi$ and plays a key role in compressed sensing \cite{bacami, doela, elabru, gilmutstr, trop04}. Moreover, for fixed $M$,  the minimizers of $\FP_{\infty, M}$   are called Grassmanian frames\index{frame!Grassmanian frames} \cite{beko, strhea}. By using continuity and compactness arguments one can show that given $M, N,$ $\FP_{\infty, M}$  always has  a minimum \cite[Appendix]{beko}. The challenge is the construction of these Grassmanian. In \cite{beko} constructions of Grassmanians were considered for $N=2$ and  $M\geq 2,$  and for $N=3$  when $M\in \{3,4,5,6\}$. The ideas used in these constructions are based on analytical interpretation of some geometric results obtained in \cite{jlt}. The general question of constructing the minimizers of $\FP_{\infty, M}$ for $N\geq 3$ and $M\geq 6$ is still a mostly open question.

Even more, minimizing $\FP_{p, M}$ is an extremely difficult problem as one needs to deal with both $p, M,$ and the ambient dimension $N$.   Some results on the minimizers as well as the value of the minimum as a function of the involved parameters were proved in \cite{ehloko10}.  We refer to \cite{Oktay} for earlier results on minimizing the $p^{th}$ frame potential. Before summarizing some of these results we consider the special case where $M=3$, $N=2$ and seek the minimizers of $$\FP_{p, 3}(\Phi)=\sum_{k, \ell =1}^3 |\ip{\varphi_k}{\varphi_\ell}|^p$$ when $p\in (0, \infty)$ with the usual modification when $p=\infty$, and $\Phi=\{\varphi_k\}_{k=1}^3 \subset S^1$. 

When $p=2$, $$\FP_{2, 3}(\Phi)=\sum_{k, \ell =1}^3 |\ip{\varphi_k}{\varphi_\ell}|^2\geq 9/2$$ with equality if and only if $\Phi=\{\varphi_k\}_{k=1}^3 \subset S^1$ is a FUNTF.  A minimizer of $\FP_{2, 3}$ is the MB-frame which is pictured below:

\begin{figure}[h]
\begin{center}
\includegraphics[scale=0.5]{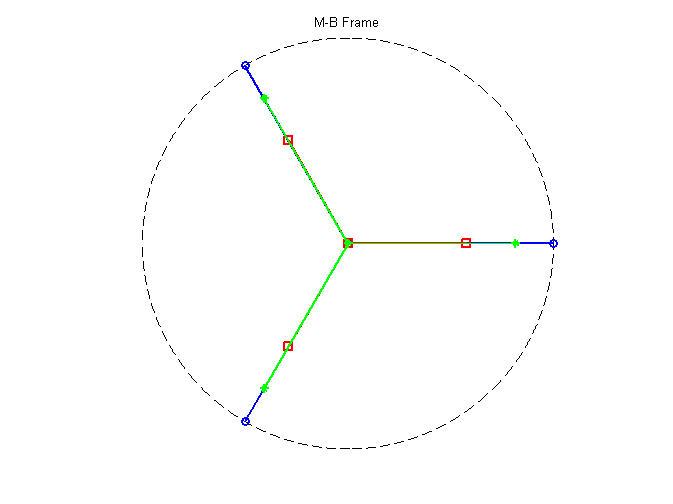}
\end{center}
\caption{An example of Equiangular FUNTF: the MB-frame.}
\label{fig-6}
\end{figure}

When $p=\infty$, $$\FP_{\infty, 3}(\Phi)=\max_{k \neq \ell}|\ip{\varphi_k}{\varphi_\ell}|\geq 1/\sqrt{2}$$ with equality if and only if $\Phi=\{\varphi_k\}_{k=1}^3 \subset S^1$ is an ETF. But what happens for other values of $2\neq p \in (0, \infty)$? This was partially answered in \cite{ehloko10} for $0<p\leq 2$ and recently the case $p\geq 2$ was settled \cite{ahamaps-reu}. Before giving more details on this case, we first collect a number of generic results about the minimizers of $\FP_{p,M}$ when $M\geq N\geq  2$ and $p \in (0, \infty]$.

\begin{proposition}\label{recapminpfp} Let $p \in (0, \infty]$, $M, N$ be  positive integers. Let $\Phi=\{\varphi_k\}_{k=1}^M \subset S^{N-1}$ we have:
\begin{enumerate}
\item[(a)] If $M\geq N$ and $2< p < \infty,$ then $$\FP_{p,M}(\Phi)\geq M(M-1)\big(\tfrac{M-N}{N(M-1)}\big)^{p/2} +N,$$  and equality holds if and only if $\Phi$ is an ETF.
\item[(b)] Let $0<p<2$ and assume that $M=kN$ for some positive integer $k$. Then the minimizers of the $p$-frame potential are exactly the $k$ copies of any orthonormal basis modulo multiplications by $\pm 1$. The minimum of~\eqref{eq:pth frame potential} over all sets of $M=kN$ unit norm vectors is $k^{2}N$.
\item[(c)] Assume that $M=N+1$ and set $p_{0}=\frac{\log(\frac{N(N+1)}{2})}{\log(N)}= \frac{\log(\frac{M(M-1)}{2})}{\log(M-1)}$. Assume that $\FP_{p_{0},N}(\Phi) \geq N+3,$ with equality holding if and only if $\Phi=\{\varphi_k\}_{k=1}^{N+1}$ is an orthonormal basis plus one repeated vector or an equiangular FUNTF. Then,
\begin{itemize}
\item[\textnormal{(1)}] for $0<p<p_{0} $,  for any $\Phi=\{\varphi_k\}_{k=1}^{N+1}\subset S^{N-1}$, we have  $\FP_{p,N+1}(\Phi) \geq N+3,$ and equality holds if and only if $\Phi=\{\varphi_k\}_{k=1}^{N+1}$ is an orthonormal basis plus one repeated vector,
\item[\textnormal{(2)}] for $p_{0}<p<2$,  for any $\Phi=\{\varphi_k\}_{k=1}^{N+1}\subset S^{N-1}$, we have $\FP_{p,N+1}(\Phi) \geq 2^{\frac{p}{p_{0}}}\, ((N+1)N)^{1-\frac{p}{p_{0}}} + N+1,$ 
and equality holds if and only if $\Phi=\{\varphi_k\}_{k=1}^{N+1}$ is an equiangular FUNTF.
\end{itemize}
\end{enumerate}
\end{proposition}

In the special case when $N=2$, part (c) of the proposition becomes:

\begin{cor}\label{cor:fpp32}
For $N=2$, $M=3$, and $p_{0}=\frac{\log(3)}{\log(2)}$,  the hypothesis of (c) above holds. That is, for any $\Phi=\{\varphi_k\}_{k=1}^3 \subset S^{1}$, $$\FP_{p_{0},3}(\Phi) \geq 5,$$ and equality holds if and only if $\Phi=\{\varphi_k\}_{k=1}^3$ is an orthonormal basis plus one repeated vector or an equiangular FUNTF. 
\end{cor}

However, when $N\geq 3,$ it is still unknown if the hypothesis of proposition~\ref{recapminpfp} (c) holds, and it was conjectured in in \cite{ehloko10} that with $p_0$ given in (c), $$\FP_{p_{0}, N+1}(\Phi)\geq N+3$$ with equality if and only if $\Phi=\{\varphi_k\}_{k=1}^{N+1}$ is an orthonormal basis plus one repeated vector or an equiangular FUNTF.

Using Corollary~\ref{cor:fpp32} one can compute   $$\mu_{p, 3, 2}=\min \{\FP_{p, 2}(\Phi): \Phi=\{\varphi_k\}_{k=1}^3 \subset S^1\}$$  for all $p \in (0, \infty]$ leading to  

$$\mu_{p, 3, 2}= \begin{cases}
5 & \textrm{for} \,  p \in (0, \frac{\log(3)}{\log(2)}]\\
3+ 6 e^{-p\log 2} & \textrm{for} \, p \geq \frac{\log(3)}{\log(2)}.
\end{cases}$$

\begin{figure}[h]
\begin{center}
\includegraphics[scale=0.7]{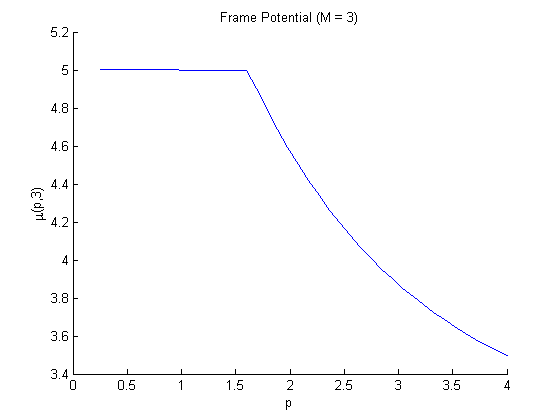}
\end{center}
\caption{Graph of $\mu_{p,3,2}$ when $p \in (0, 4)$.}
\label{fig-7}
\end{figure}

The graph of $\mu_{p,3,2}$ when $p \in (0, \infty)$ is given  in Figure~\ref{fig-7}. 

\begin{figure}[h]
\begin{center}
\includegraphics[scale=0.7]{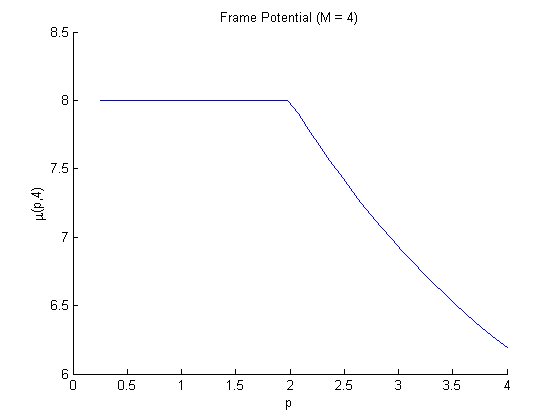}
\end{center}
\caption{Graph of $\mu_{p,4,2}$ when $p \in (0, 4)$.}
\label{fig-8}
\end{figure}

\begin{figure}[h]
\begin{center}
\includegraphics[scale=0.7]{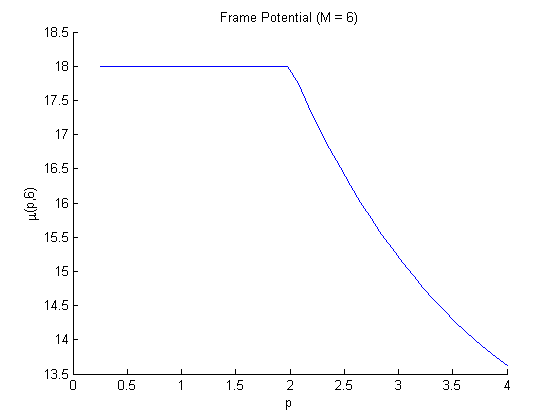}
\end{center}
\caption{Graph of $\mu_{p,6,2}$ when $p \in (0, 4)$.}
\label{fig-9}
\end{figure}

One can ask about  of $\mu_{p, M, 2}$  for other values of $M.$ It follows from proposition~\ref{recapminpfp}  (b) that $\mu_{p,M,2}=2k^2$ for all $p \in (0,2]$ whenever $M=2k$ is an even integer. For $p>2$ or odd $M$  some numerical simulations were  considered in \cite{ahamaps-reu}. For example,  the following graphs (Figures~\ref{fig-8} and~\ref{fig-9})  of  $\mu_{p, M, 2}$ for $M \in \{4,  6\}$ were obtained.

\begin{figure}[h]
\begin{center}
\includegraphics[scale=0.7]{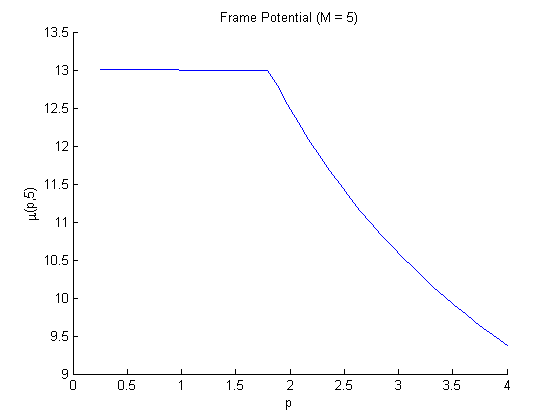}
\end{center}
\caption{Graph of $\mu_{p,5,2}$ when $p \in (0, 4)$.}
\label{fig-10}
\end{figure}

\begin{figure}[h]
\begin{center}
\includegraphics[scale=0.7]{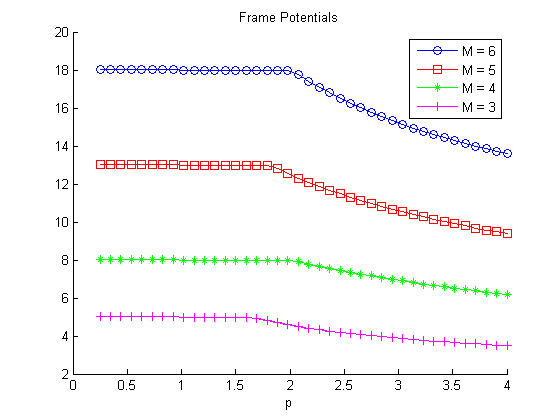}
\end{center}
\caption{Graph of $\mu_{p,M,2}$ when $p \in (0, 4)$, $M \in \{3, 4, 5, 6\}$.}
\label{fig-11}
\end{figure}

For $M=5$ the numerical results suggest that the graph of $\mu_{p, 5, 2}$ is as given in~\ref{fig-10}.  Finally, we plot the behavior of $\mu_{p, M, 2}$ as a sequence in $M$  when $p \in (0, 4)$ is shown in Figure~\ref{fig-11}. 

For integer values of $p$, the minimizers of $\FP_{p,M}$ have been investigated in connection with the theory of spherical designs \cite{Delsarte:1977aa, Venkov:2001aa}.

\begin{definition}\label{tdesign}
Let $t$ be a positive integer. 
A \emph{spherical $t$-design}\index{spherical $t$-design} is a finite subset $\{x_i\}_{i=1}^M$ of the unit sphere $S^{N-1}$ in $\R^N$,
such that,
\begin{equation*}
\frac{1}{M}\sum_{i=1}^M h(x_i) = \int_{S^{N-1}} h(x)d\sigma(x),
\end{equation*}
for all homogeneous polynomials $h$ of total degree equals or less than $t$ in $N$ variables and where $\sigma$ denotes the uniform surface measure on $S^{N-1}$ normalized to have mass one. 
\end{definition}

It is easy to see that any spherical $t-$design is also a spherical $t'-$design for all positive integers $t'\leq t$. Spherical $2-$designs are exactly FUNTFs whose center of mass is at the origin. More precisely we have: 

\begin{proposition}\label{2designfuntf}
$\Phi=\{\varphi_k\}_{k=1}^M \subset S^{N-1}$ is a spherical $2$-design if and only if $\Phi$ is a FUNTF and $\sum_{k=1}^M\varphi_k=0.$
\end{proposition}

We refer to \cite{Delsarte:1977aa},  \cite[Theorem 3.2]{Venkov:2001aa} for details on the proof of the above proposition. Recalling that FUNTFs minimize the frame potential, it is not surprising that spherical $t$-designs also minimize a potential. In particular,

 \begin{theorem}\label{theorem:p even integer discrete}\cite[Theorem 8.1]{Venkov:2001aa}
Let $p=2k$ be an even integer and $\{x_i\}_{i=1}^M=\{-x_i\}_{i=1}^M\subset S^{N-1}$, then 
  \begin{equation*}
\FP_{p, M}(\{x_{i}\}_{i=1}^{M})  \geq  \frac{1\cdot 3\cdot 5\cdots(p-1)}{N(N+2)\cdots (N+p-2) }M^2,
 \end{equation*}
and equality holds if and only if $\{x_i\}_{i=1}^M$ is a spherical $p$-design. 
 \end{theorem} 

\subsection{ Probabilistic $p$ frame potential}\label{subsec3.4}
The $p^{th}$ frame potential can be viewed in light of mass distributions on the unit sphere. It is therefore natural to look at it from a probabilistic point of view. This motivates the  the introduction of the larger family of potential called \emph{probabilistic $p$-frame potential}\index{$p$-frame potential!probabilistic $p$-frame potential}. 

For $p \in (0, \infty)$ set $$\mathcal{P}_{p}=\big\{ \mu \in \mathcal{P}: M_{p}^{p}(\mu)=\int_{\R^N}\|y\|^p d\mu(y) < \infty\big\}.$$

\begin{definition}\label{pprobpot}
For  each $p \in (0, \infty),$ the \emph{probabilistic $p-$frame potential} is given by
\begin{equation}\label{peqpfp}
\pfp(\mu, p)=\iint_{\R^{N}\times \R^{N}}|\langle x,y\rangle |^{p}\, d\mu(x)\, d\mu(y).
\end{equation}
\end{definition}
When $\mu$ is a purely atomic measure with atoms on the unit sphere, that is when $\supp(\mu)=\Phi=\{\varphi_k\}_{k=1}^{M} \subset S^{N-1},$  $\pfp(\mu, p)$ reduces to $\FP_{p,M}$ given in~\eqref{eq:pth frame potential}. 

This class of potentials is related to  the potentials considered by G.~Bj\"rock   \cite{Bjoerck:1955aa}. More precisely, suppose $F\subset \R^N$ is compact and let $\lambda>0$.  Bj\"orck considered the question of maximizing the functional  $$I_F(\mu)=\iint_{F\times F}\|x-y\|^{\lambda}d\mu(x)d\mu(y)$$ where  $\mu$ ranges over all  positive Borel measures with $\mu(F)=1$. It turns out that the techniques used in  \cite{Bjoerck:1955aa} to maximize  $I_F(\mu)$ can be extended to understand the minimizers of  $\pfp(\mu)$ when $\mu$ is restricted to a probability measure on the unit sphere $S^{N-1}$ in $\R^N$. In particular,  it was proved in \cite[theorem 4.9]{ehloko10} that when restricted to probability measures $\mu$ supported on the unit sphere of $\R^N$  and when $0<p<2$, then the minimizers of  $\pfp(\mu, p)$ are discrete probability measures. Furthermore, the support of such minimizers contains an orthonormal basis $B$ and is contained in the set $\pm B$. More specifically we have:

 \begin{theorem}\cite[Theorem 4.9]{ehloko10}\label{theorem:p less than 2}
 Let $0<p<2$, then the minimizers of  \eqref{peqpfp} over all the probability measures supported on the unit sphere $S^{N-1}$ are exactly those probability measures $\mu$ that satisfy
 \begin{itemize}
 \item[(i)] there is an orthonormal basis $\{e_1,\ldots,e_N\}$ for $\R^N$ such that 
 \begin{equation*}
 \{e_1,\ldots,e_N\} \subset \supp(\mu) \subset \{\pm e_1,\ldots,\pm e_N\}
 \end{equation*}
 \item[(ii)] there is $f:S^{N-1}\rightarrow \R$ such that $\mu(x) = f(x) \nu_{\pm x_1,\ldots,\pm x_N}(x)$ and
 \begin{equation*}
 f(x_i) +f(-x_i) = \frac{1}{N},
 \end{equation*} 
where the measure $\nu_{\pm x_1,\ldots,\pm x_N}(x)$ represent  the counting measure of the set $\{\pm x_i : i=1,\ldots,N\}$.
 \end{itemize} 
 \end{theorem}

Theorem~\ref{theorem:p even integer discrete} shows that the minimizers of $\FP_{p, M}$ when $p=2k$ is an even integer, are exactly spherical $p-$designs. In view of this fact, one can ask whether the minimizers of $\pfp$ have some special ``approximation'' properties. This partially motivates that following definition in which   we denote by $\mathcal{M}(S^{N-1},\mathcal{B})$ the space of all Borel probability measures supported on $S^{N-1}$.

\begin{definition}\cite[Definition 4.1]{ehloko10}\label{probpframe}
For $0<p<\infty$, we call $\mu\in\mathcal{M}(S^{N-1},\mathcal{B})$ a \emph{probabilistic $p$-frame}\index{probabilistic $p$-frame} for $\R^N$ if and only if there are constants $A,B>0$ such that
 \begin{equation}\label{ppframeineq}
A\|y\|^p\leq  \int_{S^{N-1}} |\ip{ x}{y}|^p d\mu(x) \leq B\|y\|^p, \quad\forall y\in\R^N.
 \end{equation}
 We call $\mu$ a \emph{tight probabilistic $p$-frame}\index{probabilistic $p$-frame!tight probabilistic $p$-frame} if and only if we can choose $A=B$.
 \end{definition}
By symmetry considerations, it is not difficult to show that  the uniform surface measure $\sigma$ on $S^{N-1}$ is always a tight probabilistic $p$-frame, for each $0<p<\infty$.  In addition, observe that  we can always take $B=1$ in~\eqref{ppframeineq}. Thus to determine if a probability measure $\mu$ on $S^{N-1}$ is a probabilistic $p-$frame one must focus on establishing the lower bound in the above definition.   When $p=2$ this definition reduces to that of  probabilistic frames introduced earlier. In fact, more is true: 
 
 \begin{lemma}\cite[Lemma 4.5]{ehloko10}\label{equivappf}  If $\mu$ is probabilistic frame, then it is a probabilistic $p$-frame for all $1\leq p< \infty$. Conversely, if $\mu$ is a probabilistic $p$-frame for some $1\leq p < \infty$, then it is a probabilistic frame.
\end{lemma}

The analogy between tight probabilistic $p$-frames and spherical $t-$designs can now be made explicitly as one can show the following result which is an analog of theorem~\ref{theorem:p even integer discrete}. More specifically, the result below shows that  tight probabilistic $p$-frames are the minimizers of the probabilistic frame potential~\eqref{peqpfp} when restricted to probability measures supported on $S^{N-1}$, and when $p$ is an even integer:

\begin{theorem}\cite[Theorem 4.10]{ehloko10}\label{theorem:p even integer}
 Let $p$ be an even integer. For any probability measure $\mu$  on $S^{N-1}$,  
  \begin{equation*}
 \PFP(\mu, p)=\int_{S^{N-1}}\int_{S^{N-1}} |\ip{x}{y}|^p d\mu(x) d\mu(y) \geq  \frac{1\cdot 3\cdot 5\cdots(p-1)}{N(N+2)\cdots (N+p-2) },
 \end{equation*}
and equality holds if and only if $\mu$ is a tight probabilistic  $p$-frame. 
 \end{theorem}

 By combining Theorem~\ref{theorem:p even integer discrete} and Theorem~\ref{theorem:p even integer}  we can conclude that when $p=2k$ there exists a one-to-one correspondence between the class of spherical $p-$designs and the class of discrete tight probabilistic $p-$frames. More specifically,  every spherical $p-$design supports  a discrete measure $\mu$ which is a tight probabilistic $p$-frame. This is summarized in the following proposition:
 
 \begin{proposition} Let $p=2k$ be an even positive integer. A set $\Phi=\{\varphi_k\}_{k=1}^M\subset S^{N-1}$ is a spherical $p-$design if and only if the probability measure $\mu_\Phi=\tfrac{1}{M}\sum_{k=1}^N\delta_{\varphi_{k}}$ is a tight probabilistic $p-$frame. 
 
 \end{proposition}

 The question then becomes how to construct  tight probabilistic $p$-frames. When restricted to discrete measures and when $p=2k$ is an even integer, this problem  is equivalent to constructing spherical $p-$designs. This is a difficult problem with known solutions only for certain values of $p, M, $ and $N$. Of course, and as shown in Section~\ref{sec2} the special case when $p=2$ leads to the  FUNTFs.  
 The analytics methods developed in \cite{cw14} are new promising techniques that could be used to investigated in general the minimizers of $\pfp(\mu)$ when $\mu$ ranges over the probability measures on $S^{N-1}$ and $p>0$.

\section*{Acknowledgment}
This work was partially supported by a grant from the Simons Foundation ($\# 319197$ to Kasso Okoudjou).
The author would like to  thank Chae Clark and Matthew Begu\'e for their helpful discussions.

\end{document}